\documentclass[12pt]{article}

\usepackage{fullpage,graphicx}
\usepackage{amsmath, amssymb,amsthm}
\usepackage[hidelinks]{hyperref}
\usepackage{tikz}
\usetikzlibrary{positioning,calc}
\usetikzlibrary{automata,positioning,shapes,shapes.misc,shapes.geometric}
\usetikzlibrary{math}
\usetikzlibrary {arrows.meta}
\usepackage{caption}
\usepackage{subcaption}

\newtheorem{thm}{Theorem}[section]
\newtheorem{lem}[thm]{Lemma}
\newtheorem{cor}[thm]{Corollary}

\newtheorem{example}[thm]{Example}
\newtheorem{remark}[thm]{Remark}
\theoremstyle{definition}
\newtheorem{defn}[thm]{Definition}

\newcommand{\ints}{\mathbb Z}

\newcommand{\re}{\mathbb R}

\newcommand{\cK}{\mathcal K}

\newcommand{\G}{\Gamma}

\begin{document}

\title{On Lewis's homotopy in distance-regular graphs}
\author{MollyAnn Burkey, William J.\ Martin, Freud Oulon\\
Department of Mathematical Sciences, Worcester Polytechnic Institute \\
100 Institute Road, Worcester, Massachusetts, USA 01609}
\date{\today}

\maketitle

\begin{abstract}
In 2000, Lewis \cite{lewis} introduced a nested sequence of normal subgroups of the fundamental
group of a distance-regular graph $\Gamma$. While the full homotopy group itself is simply a free group, it is
interesting to determine the smallest integer $k$ such that every closed walk in $\Gamma$ can be built, up
to homotopy, from closed walks of ``essential length'' $k$ or less. Lewis proved that $k\le 6$ for any
$Q$-polynomial distance-regular graph $\Gamma$ that is not a pseudoquotient. In this paper, we prove
that $k=3$ for Johnson graphs and that $k=4$ in the case of Hamming graphs. We also compute the homotopy
groups of natural cell  complexes associated to the folded Johnson graphs and folded cubes. The paper
finishes with some open questions and suggestions for further lines of investigation.
\end{abstract}

\section{Motivation and background}
\label{Sec:intro}

We begin with our motivation for this investigation. We then  briefly summarize background material and review Lewis's notation and results.

Distance-regular graphs \cite{del,banito,bcn, godsil, vDKT} are fundamental in algebraic graph theory, with rich connections
to group actions, extremal graph theory, coding theory, combinatorial design theory and more.  $Q$-polynomial distance-regular
graphs have particularly rich structure, including some with Terwilliger algebras tied to the representation theory of Lie algebras
and quantum groups \cite{itotanter}. A number of families of
distance-regular graphs admit some sort of local characterization. For example, any graph $\G$  which is locally a grid
and has $| \G(x) \cap \G(y) | =4$ for any vertices $x,y$ at distance two is either a Johnson graph or a folded Johnson
graph. In \cite{lewis}, Lewis introduced another local-to-global notion. While the fundamental group $\pi(\G,x)$ of a graph
$\G$ with respect to basepoint $x$ is not particularly interesting in  itself --- it is merely a free group on  $vk/2+1-v$
generators for a $k$-regular graph on $v$ vertices ---  there is a nested sequence  of normal subgroups $G_i=\pi(\G,x,i)$
  that tell us the extent to which local topological information determines global information; Essential $G_i$ is the subgroup of $\pi(\G,x)$ generated by equivalence classes containing a representative of (essential) length $i$ or less. The corresponding factor groups $\pi(\G,x)/G_i$ can be viewed as the fundamental groups of cell complexes obtained by gluing 2-cells across all cycles in $\G$ of length
at most $i$. An interesting problem is to determine which of these cell complexes are simply connected.
(An alternative approach is taken in \cite{shult,SuzukiHomotopy} where universal covers are a core objective.)

Given a distance-regular graph $\G$, let us denote its girth by $g_1(\G)$. It is not hard to see that $G_i=1$ for $i<g_1(\G)$.  Let $g_2(\G)$ be the smallest integer $i$ for which $G_i = \pi(\G,x)$. Obviously $3\le g_1(\G)\le g_2(\G)$. In \cite[Theorem 26]{lewis},  Lewis proved that, for any  connected graph $\G$ of diameter $d$, $g_2(\G)\le 2d+1$. Unpublished conjectures of Suzuki and Koolen assert that any distance-regular graph of valency at least three has girth at most 12 (See Problem 32 in \cite[Section 18.7]{vDKT}.)  While we expect $g_2(\G)\le 12$ for
most distance-regular graphs, this is not always the case. Among the known distance-regular graphs  for which $\pi(\G,x,12) \neq \pi(\G,x)$, we explore the folded cubes and folded Johnson graphs in
this paper.

Lewis introduced the concept of a \emph{pseudoquotient} (see Definition \ref{Def:pseuodoquotient} below) and
proved that, in any $Q$-polynomial distance-regular graph of valency at least three that is not a pseudoquotient, every closed
walk can be decomposed, into closed walks of essential length six or less: that is, $\pi(\G,x)=\pi(\G,x,i)$.

Earlier, Brouwer, Cohen and Neumaier
had proved that the numerical girth of a $Q$-polynomial distance-regular graph of valency at least three is at most seven
\cite[Theorem 8.3.6]{bcn} and conjectured that the girth of such a graph should be at most six, with equality only for the Odd graphs.
In \cite{lewis}, Lewis also proved that any $Q$-polynomial distance-regular graph of valency at least three has
girth at most six \cite[Corollary 30]{lewis} by leveraging Terwilliger's balanced set condition.
Recently, Miklavi\v{c} resolved the second part of the conjecture of Brouwer, et al.\ by determining
those $Q$-polynomial distance-regular graphs having girth exactly six: they are the Odd graphs and point graphs of thin generalized hexagons (i.e., incidence graphs of finite projective planes).

This prompts us to explore $G_i=\pi(\G,x,i)$ and $\pi(\G,x)/G_i$, $i=3,\ldots,6$, for specific families of $Q$-polynomial distance-regular graphs. In this paper, we determine these homotopy groups and factor groups for Johnson graphs, Hamming graphs, folded Johnson graphs and folded hypercubes.

\section{Lewis's formulation of homotopy on graphs}
\label{Subsec:lewis}
The purpose of this section is to summarize the material in \cite{lewis} essential to our goal. The reader may consult \cite{lewis} for further details.

Let us begin with an arbitrary simple graph $\G=(X,R)$ with vertex set $X$ and adjacency relation $R$, an irreflexive and symmetric binary relation on $X$. For $x,y\in X$ and $\ell \ge 0$, a \emph{walk} in $\G$
of length $\ell$ with \emph{initial vertex} $x$ and \emph{terminal vertex} $y$ is a sequence
\[
 p = (x_0,x_1,\ldots, x_\ell)
\]
of vertices in $\G$ such that $x_0=x$, $x_\ell=y$ and $(x_{i-1},x_i)\in R$ for each $1\le i\le \ell$. Walk $p$ is
a \emph{closed walk} if $x_0=x_\ell$.  The \emph{reverse walk} of $p$ is the walk
\[
 p^{-1} = (x_\ell,x_{\ell-1},\ldots, x_0).
\]

Choose a \emph{base vertex} or \emph{basepoint} $x\in X$ and
let $\psi(\G,x)$ denote the set of all closed walks in $\G$ with initial and terminal vertex $x$.
The following definitions from Lewis \cite[Definitions 1--5]{lewis} give us what will be the elements our fundamental group.

\
\begin{itemize}
\item Let $p,q\in \psi(\G,x)$ with $p=(x=x_0,x_1,\dots,x_{\ell}=x)$. We say $q$ extends $p$ if there is an integer $i$ and a vertex $z\in X$ such that  $q = (x=x_0,x_1,\dots,x_{i-1},x_i,z,x_i,x_{i+1},\dots,x_\ell=x)$.
\item  Let $S = (p_0,p_1,\dots,p_n)$ be a sequence of walks in $\psi(\G,x)$. We say that $S$ is an \emph{admissible sequence} if, for each $i$, $1\le i\le n$,  $p_{i-1}$ extends  or is extended by $p_i$.
\item  Define a \emph{homotopy relation} $\sim$ on $\psi(\G,x)$ as follows:  for all $p,q\in \psi(\G,x)$, write $p\sim q$ whenever there exists a nonnegative integer $n$ and walks $p=p_0,p_1,\dots,p_n=q$ in $\psi(\G,x)$ such that the sequence $(p_0,p_1,\dots,p_n)$ is admissible. We say that $p$ and $q$ are \emph{homotopic} whenever $p\sim q$. Note that $\sim$ is an equivalence relation on $\psi(\G,x)$.
\item  For $p \in \psi(\G,x)$, let $[p]$ denote the equivalence class of $\sim$ containing $p$. Let $\pi(\G,x)$
denote the collection of all equivalence classes of $\sim$ on $\psi(\G,x)$.
For an equivalence class $u\in \pi(\G,x)$ and $p\in \psi(\G,x)$, we say $p$ is a \emph{representative}
of $u$ if $u = [p]$.
\end{itemize}

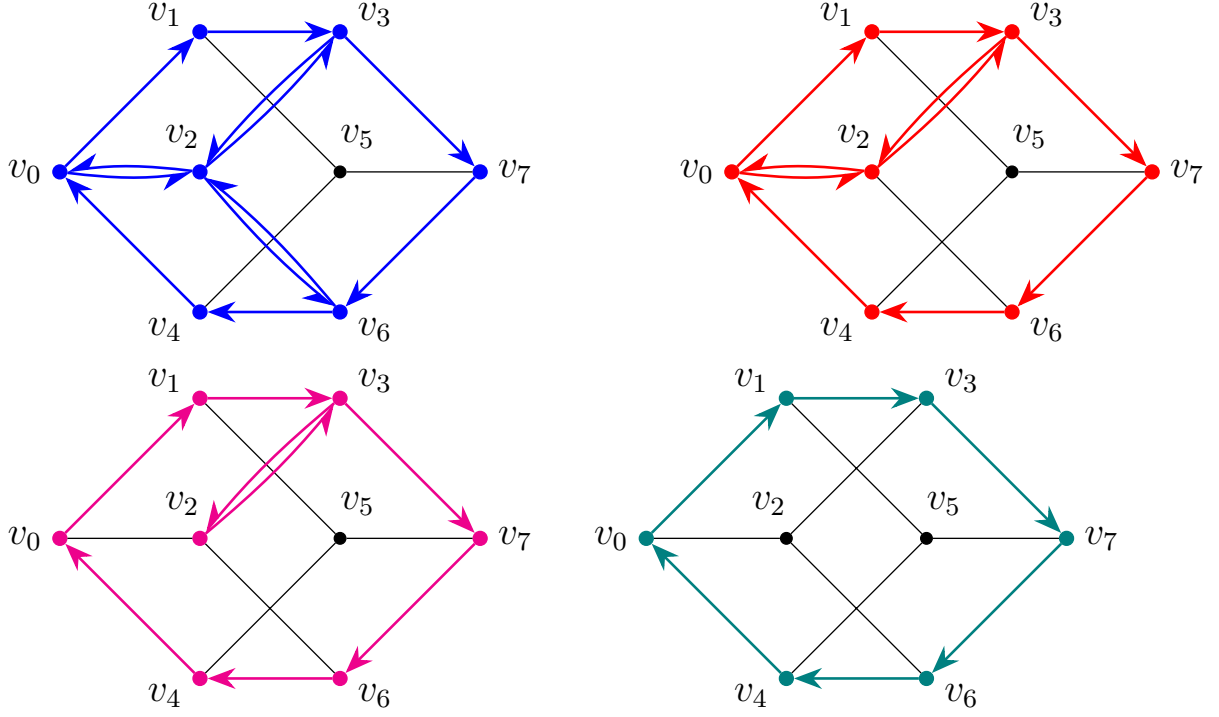
\begin{figure} [h]
\resizebox{3in}{!}{
\begin{tikzpicture}[xscale=0.75,yscale=0.75,
    vert/.style={draw,fill=black,circle,inner sep=1.5pt},
    black_vert/.style={draw,fill=black,circle,inner sep=1.25pt}]]

    \node[vert, blue] (v0) at (-3,0) {}; \node (lv0) at (-3.5,0) {$v_0$};
    \node[vert, blue] (v4) at (-1,-2) {}; \node (lv4) at (-1.5,-2.25) {$v_4$};
    \node[vert, blue] (v2) at (-1,0) {}; \node (lv2) at (-1.25,0.5) {$v_2$};
    \node[vert, blue] (v6) at (1,-2) {}; \node (lv6) at (1.5,-2.25) {$v_6$};
    \node[vert, blue] (v1) at (-1,2) {}; \node (lv1) at (-1.5,2.25) {$v_1$};
    \node[black_vert] (v5) at (1,0) {}; \node (lv5) at (1.25,0.5) {$v_5$};
    \node[vert, blue] (v3) at (1,2) {}; \node (lv3) at (1.5,2.25) {$v_3$};
    \node[vert, blue] (v7) at (3,0) {}; \node (lv7) at (3.5,0) {$v_7$};

    \draw
        (v5) -- (v1)
        (v5) -- (v7)
        (v5) -- (v4);

    \draw [thick, blue, -{Stealth[length=3mm]}]
    (v0) -- (v1);
    \draw [thick, blue, -{Stealth[length=3mm]}]
    (v1) -- (v3);
    \draw[blue, thick, bend right=5, -{Stealth[width=3mm,length=3mm,harpoon,swap]}]
        (v3) to node [auto] {} (v2);
    \draw[blue, thick, bend right=7, -{Stealth[width=3mm,length=3mm,harpoon,swap]}]
        (v2) to node [auto] {} (v0);
    \draw[blue, thick, bend right=7, -{Stealth[width=3mm,length=3mm,harpoon,swap]}]
        (v0) to node [auto] {} (v2);
    \draw[blue, thick, bend right=5, -{Stealth[width=3mm,length=3mm,harpoon,swap]}]
        (v2) to node [auto] {} (v3);
    \draw [thick, blue, -{Stealth[length=3mm]}]
        (v3) -- (v7);
    \draw [thick, blue, -{Stealth[length=3mm]}]
        (v7) -- (v6);
    \draw[blue, thick, bend right=5, -{Stealth[width=3mm,length=3mm,harpoon,swap]}]
        (v6) to node [auto] {} (v2);
    \draw[blue, thick, bend right=5, -{Stealth[width=3mm,length=3mm,harpoon,swap]}]
        (v2) to node [auto] {} (v6);
    \draw [thick, blue, -{Stealth[length=3mm]}]
        (v6) -- (v4);
    \draw [thick, blue, -{Stealth[length=3mm]}]
        (v4) -- (v0);

\end{tikzpicture}
}
\resizebox{3in}{!}{
\begin{tikzpicture}[xscale=0.75,yscale=0.75,
    vert/.style={draw,fill=black,circle,inner sep=1.5pt},
    black_vert/.style={draw,fill=black,circle,inner sep=1.25pt}]]

    \node[vert, red] (v0) at (-3,0) {}; \node (lv0) at (-3.5,0) {$v_0$};
    \node[vert, red] (v4) at (-1,-2) {}; \node (lv4) at (-1.5,-2.25) {$v_4$};
    \node[vert,red ] (v2) at (-1,0) {}; \node (lv2) at (-1.25,0.5) {$v_2$};
    \node[vert, red] (v6) at (1,-2) {}; \node (lv6) at (1.5,-2.25) {$v_6$};
    \node[vert, red] (v1) at (-1,2) {}; \node (lv1) at (-1.5,2.25) {$v_1$};
    \node[black_vert] (v5) at (1,0) {}; \node (lv5) at (1.25,0.5) {$v_5$};
    \node[vert, red] (v3) at (1,2) {}; \node (lv3) at (1.5,2.25) {$v_3$};
    \node[vert, red] (v7) at (3,0) {}; \node (lv7) at (3.5,0) {$v_7$};

    \draw
        (v5) -- (v4)
        (v5) -- (v7)
        (v5) -- (v1)
        (v2) -- (v6);

    \draw [thick, red, -{Stealth[length=3mm]}]
        (v0) -- (v1);
    \draw [thick, red, -{Stealth[length=3mm]}]
        (v1) -- (v3);
    \draw[red, thick, bend right=5, -{Stealth[width=3mm,length=3mm,harpoon,swap]}]
        (v3) to node [auto] {} (v2);
    \draw[red, thick, bend right=7, -{Stealth[width=3mm,length=3mm,harpoon,swap]}]
        (v2) to node [auto] {} (v0);
    \draw[red, thick, bend right=7, -{Stealth[width=3mm,length=3mm,harpoon,swap]}]
        (v0) to node [auto] {} (v2);
    \draw[red, thick, bend right=5, -{Stealth[width=3mm,length=3mm,harpoon,swap]}]
        (v2) to node [auto] {} (v3);
    \draw [thick, red, -{Stealth[length=3mm]}]
        (v3) -- (v7);
    \draw [thick, red, -{Stealth[length=3mm]}]
        (v7) -- (v6);
    \draw [thick, red, -{Stealth[length=3mm]}]
        (v6) -- (v4);
    \draw [thick, red, -{Stealth[length=3mm]}]
        (v4) -- (v0);

\end{tikzpicture}
}
\resizebox{3in}{!}{
\begin{tikzpicture}[xscale=0.75,yscale=0.75,
    vert/.style={draw,fill=black,circle,inner sep=1.5pt},
    black_vert/.style={draw,fill=black,circle,inner sep=1.25pt}]]

    \node[vert, magenta] (v0) at (-3,0) {}; \node (lv0) at (-3.5,0) {$v_0$};
    \node[vert, magenta] (v4) at (-1,-2) {}; \node (lv4) at (-1.5,-2.25) {$v_4$};
    \node[vert, magenta] (v2) at (-1,0) {}; \node (lv2) at (-1.25,0.5) {$v_2$};
    \node[vert, magenta] (v6) at (1,-2) {}; \node (lv6) at (1.5,-2.25) {$v_6$};
    \node[vert, magenta] (v1) at (-1,2) {}; \node (lv1) at (-1.5,2.25) {$v_1$};
    \node[black_vert] (v5) at (1,0) {}; \node (lv5) at (1.25,0.5) {$v_5$};
    \node[vert, magenta] (v3) at (1,2) {}; \node (lv3) at (1.5,2.25) {$v_3$};
    \node[vert, magenta] (v7) at (3,0) {}; \node (lv7) at (3.5,0) {$v_7$};

    \draw
        (v5) -- (v4)
        (v5) -- (v7)
        (v5) -- (v1)
        (v2) -- (v6)
        (v0) -- (v2);

    \draw [thick, magenta, -{Stealth[length=3mm]}]
        (v0) -- (v1);
    \draw [thick, magenta, -{Stealth[length=3mm]}]
        (v1) -- (v3);
    \draw[magenta, thick, bend right=5, -{Stealth[width=3mm,length=3mm,harpoon,swap]}]
        (v3) to node [auto] {} (v2);
    \draw[magenta, thick, bend right=5, -{Stealth[width=3mm,length=3mm,harpoon,swap]}]
        (v2) to node [auto] {} (v3);
    \draw [thick, magenta, -{Stealth[length=3mm]}]
        (v3) -- (v7);
    \draw [thick, magenta, -{Stealth[length=3mm]}]
        (v7) -- (v6);
    \draw [thick, magenta, -{Stealth[length=3mm]}]
        (v6) -- (v4);
    \draw [thick, magenta, -{Stealth[length=3mm]}]
        (v4) -- (v0);

\end{tikzpicture}
}
\resizebox{3in}{!}{
\begin{tikzpicture}[xscale=0.75,yscale=0.75,
    vert/.style={draw,fill=black,circle,inner sep=1.5pt},
    black_vert/.style={draw,fill=black,circle,inner sep=1.25pt}]]

    \node[vert, teal] (v0) at (-3,0) {}; \node (lv0) at (-3.5,0) {$v_0$};
    \node[vert, teal] (v4) at (-1,-2) {}; \node (lv4) at (-1.5,-2.25) {$v_4$};
    \node[black_vert] (v2) at (-1,0) {}; \node (lv2) at (-1.25,0.5) {$v_2$};
    \node[vert, teal] (v6) at (1,-2) {}; \node (lv6) at (1.5,-2.25) {$v_6$};
    \node[vert, teal] (v1) at (-1,2) {}; \node (lv1) at (-1.5,2.25) {$v_1$};
    \node[black_vert] (v5) at (1,0) {}; \node (lv5) at (1.25,0.5) {$v_5$};
    \node[vert, teal] (v3) at (1,2) {}; \node (lv3) at (1.5,2.25) {$v_3$};
    \node[vert, teal] (v7) at (3,0) {}; \node (lv7) at (3.5,0) {$v_7$};

    \draw
        (v5) -- (v4)
        (v5) -- (v7)
        (v5) -- (v1)
        (v2) -- (v6)
        (v0) -- (v2)
        (v2) -- (v3);

    \draw [thick, teal, -{Stealth[length=3mm]}]
        (v0) -- (v1);
    \draw [thick, teal, -{Stealth[length=3mm]}]
        (v1) -- (v3);
    \draw [thick, teal, -{Stealth[length=3mm]}]
        (v3) -- (v7);
    \draw [thick, teal, -{Stealth[length=3mm]}]
        (v7) -- (v6);
    \draw [thick, teal, -{Stealth[length=3mm]}]
        (v6) -- (v4);
    \draw [thick, teal, -{Stealth[length=3mm]}]
        (v4) -- (v0);

\end{tikzpicture}
}
\caption{\label{Fig:admissible sequence} {An admissible sequence of walks.}}
\end{figure}

\begin{example}
\label{Examp:3cube}
Let $\G=H(3,2)$ be the 3-cube, a $3$-regular graph on vertex set $v_0,v_1,\ldots,v_7$ where $v_j$ has as its
label the binary $3$-tuple representing the integer $j$ in base two and where two vertices are adjacent if their
labels differ in exactly one coordinate. With $x=v_0$ as basepoint and
\begin{eqnarray*}
 p_0 &=& (x=v_0,v_1,v_3,v_2,v_0,
v_2,v_3,v_7,v_6,v_2,
v_6,v_4,v_0=x) \\
 p_1 &=& (x=v_0,v_1,v_3,v_2,v_0,
v_2,v_3,v_7,v_6,v_4,v_0=x) \\
 p_2 &=& (x=v_0,v_1,v_3,v_2,v_3,v_7,v_6,v_4,v_0=x) \\
 p_3 &=& (x=v_0,v_1,v_3,v_7,v_6,v_4,v_0=x)
\end{eqnarray*}
the sequence $S=(p_0,p_1,p_2,p_3)$ is an admissible sequence of walks in $\psi( H(3,2), v_0)$.
So $p_0$ is homotopic to $p_3$, $p_0 \sim p_3$, and $[p_0]=[p_3]$. These walks are depicted in Figure \ref{Fig:admissible sequence}.
\end{example}

Lewis \cite[Lemma 6]{lewis} proves that, for any $p,q,r \in \psi(\G,x)$ such that $p$ extends $q$, $p$ extends $r$, and $q\ne r$, there exists a walk $s\in \psi(\G,x)$ such that both $q$ and $r$ extend $s$.

A closed walk $p$ with basepoint $x$ is \textit{reduced} if $p$ does not extend $q$ for any $q\in \psi(\G,x)$.

\begin{lem}[{\cite[Theorem 8]{lewis}}]
Let $\G = (X,R)$ be a graph. Fix $x\in X$, and pick any equivalence class $u \in \pi(\G,x)$. Then $u$
has exactly one reduced representative, denoted $\tilde{u}$.  Furthermore, $\tilde{u}$ is the unique representative of $u$ of minimal length.
\end{lem}

\begin{defn}[{\cite[Definition 9]{lewis}}]
Let $p=(x_0,x_1,\ldots,x_\ell)$ and $q=(y_0,y_1,\ldots,y_m)$ be walks in graph $\G$ and assume
that $x_\ell=y_0$. The \emph{concatenation} of $p$ and $q$ is the walk
$pq = (x_0,x_1,\ldots,x_\ell,y_1,\ldots,y_m)$. Observe that $pq$ is a walk in $\G$. Whenever we write $pq$, we will assume that the terminal vertex of $p$ is the initial vertex of $q$.
\end{defn}

\begin{lem}[{\cite[Lemma 10]{lewis}}]
\label{Lem:welldef}
Concatenation and the reverse operation are well-defined on elements of $\pi(\G,x)$.  Let $p,p',q,q'\in \psi(\G,x)$.
\begin{itemize}
\item[(i)] If $p\sim p'$ and $q\sim q'$, then $pq \sim p'q'$;
\item[(ii)] If $p\sim q$, then  $p^{-1} \sim q^{-1}$;
\item[(iii)] If $q\sim (x)$, then $pq\sim p$.
\end{itemize}
\end{lem}

\begin{defn}[{\cite[Definition 11]{lewis}}]
\label{Def:gpop}
Let $\G = (X,R)$ and $x\in X$ be given. For $u$ and $v$ equivalence classes in $\pi(\G,x)$, define a
binary operation (denoted simply by juxtaposition) $uv = [pq]$ where $p$ is any representative of $u$
and $q$ is any representative of $v$. By Lemma \ref{Lem:welldef}, this is a well-defined operation on
$\pi(\G,x)$.
\end{defn}

\begin{thm}[{\cite[Theorem 12]{lewis}}]
\label{Thm:fundgp}
The set $\pi(\G,x)$ forms a group under the operation given in Definition \ref{Def:gpop}:
\begin{itemize}
\item[(i)] for $u,v,w \in \pi(\G,x)$, $(uv)w=u(vw)$;
\item[(ii)] for $e=[(x)]$ we have $ue=u=eu$ for every $u \in \pi(\G,x)$;
\item[(iii)] for $u\in \pi(\G,x)$, there is a unique $u^{-1} \in \pi(\G,x)$ satisfying $u u^{-1} = e = u^{-1} u$,
namely $u^{-1}=\left[ p^{-1} \right]$ for any  representative $p$ of $u$.
\end{itemize}
\end{thm}

We will prove in Theorem $\ref{Thm:freegp}$ that $\pi(\Gamma,x)$ is a free group, but this is well known. \\
Recall that a closed walk $p$ in a graph $\G$ with initial and terminal vertex $y$  is reduced if it does not extend any walk $q\in \psi(\G,y)$. A walk
\[
p = (y=y_0,y_1,\ldots,y_{\ell-1},y_\ell=y)
\]
in $\psi(\G,y)$  is  \emph{cyclically reduced} if $\ell=0$ or $p$ is reduced with $y_1 \neq y_{\ell-1}$. Now
consider $\G=(X,R)$ with basepoint $x\in X$. Lewis proves \cite[Lemma 14]{lewis} that, for any $p\in \psi(\G,x)$,
there is a unique vertex $y$ along $p$, a unique cyclically reduced $q\in \psi(\G,y)$ and a unique walk $r$ in $\G$ such that $p=rqr^{-1}$.

\begin{defn}[{\cite[Definitions 13 and 15]{lewis}}]
\label{Def:essential}
With $\G$ and $x$ as above, pick any $u\in \pi(\G,x)$. Write $\tilde{u}=rqr^{-1}$ where $q$ is cyclically reduced.
Then the \emph{essential length} of $u$ is defined to be the length of closed walk $q$. The subgroup of
$\pi(\G,x)$ generated by all elements of essential length at most $i$ is denoted $\pi(\G,x,i)$
\end{defn}

\begin{lem}[{\cite[Lemma 18, Corollary 22, Theorem 26]{lewis}}]
\label{Lem:normal}
For each $i\ge 0$, $\pi(\G,x,i)$ is a normal subgroup of $\pi(\G,x)$. These satisfy:
\begin{itemize}
\item[(i)] $\pi(\G,x,i) \subseteq \pi(\G,x,i+1)$ for each $i\ge 0$;
\item[(ii)] $\pi(\G,x,0) = \pi(\G,x,1) = \pi(\G,x,2) = \{e\}$.
\end{itemize}
\end{lem}

The following is the main theorem of \cite{lewis}.

\begin{thm}[{\cite[Theorem 37]{lewis}}]
\label{Thm:Lewis}
Let $\G$ be a  $Q$-polynomial distance-regular graph having diameter $d\ge 3$.
\begin{itemize}
\item[(i)] $\pi(\G,x,6) \neq \{e\}$;
\item[(ii)] If $\pi(\G,x,6) \neq \pi(\G,x)$, then $\G$ is a pseudoquotient (cf.\ Definition \ref{Def:pseuodoquotient} below). In this case,
there is a unique \emph{covering diameter} $D\in \{2d,2d+1\}$  specified in Definition \ref{Def:pseuodoquotient} for which
$\pi(\G,x,6) =  \pi(\G,x,D-1) \neq \pi(\G,x,D) = \pi(\G,x)$.
\end{itemize}
\end{thm}

\section{Topology, cell complexes and quotient groups}
\label{Subsec:ess}

As a $1$-dimensional cell complex, a graph $\G$ can be viewed as a Hausdorff space and, indeed, the fundamental group of Lewis as defined above is isomorphic to the ordinary fundamental group of that Hausdorff space \cite{stallings}.

The edge set of $\G$, denoted $E(\G)$ here, is a collection of unordered pairs $\{x,y\}$ of distinct vertices in $\G$. We can also encode this as the adjacency relation $R$ via $(x,y),(y,x)\in R$ whenever $\{x,y\} \in E(\G)$.

\begin{thm}
\label{Thm:freegp}
Let $\G$ be a connected graph with $v$ vertices, $m$ edges, and let $N=m+1-v$. Then $\G$ is homotopy equivalent
to a {bouquet of circles}, a wedge sum $\displaystyle{\bigvee_{i=1}^N} S^1$ of $N$ copies of $S^1$ all joined at a common basepoint. The fundamental group $\pi_1(\G,x)$ is isomorphic to a free group on $N$ generators.
\end{thm}

\begin{proof}
Let $T$ be a spanning tree of $\G$. Since $T$ has $v-1$ edges, the number of edges of $\G$ not in $T$ is
$
m-(v-1)=m+1-v=N.
$
We will show that these $N$ edges give the $N$ circles in the bouquet.

View $\G$ as a $1$-dimensional cell complex. Since $T$ is a tree, $T$ is contractible. Collapse all of $T$ to
a single point. Let
$
q:\G \to \G/T
$
be the quotient map. Each edge of $T$ is mapped to the basepoint in the quotient, while each edge
$e\in E(\G)\setminus E(T)$ becomes a circle attached at that same point. Therefore
$
\G/T \cong \bigvee_{i=1}^N S^1.
$

The quotient map $q$ is a homotopy equivalence. One way to see this is to collapse the tree $T$ one edge
at a time, starting from leaves. Collapsing a leaf edge in a tree does not change the homotopy type of the
graph. Repeating this process collapses the whole spanning tree to a point, so $\G$ and $\G/T$ have the
same homotopy type.

It remains only to identify the fundamental group of the quotient. By example 1.21 in Hatcher \cite[page 43]{hatcher}, we have
$
\pi_1\left(\bigvee_{i=1}^N S^1,x\right) \cong F_N,
$
the free group on $N$ generators. Hence
$
\pi_1(\G,x)\cong F_N.
$
\end{proof}

Recall the concept of a \emph{cell complex}; see, e.g., Hatcher \cite[page 5]{hatcher}. In this paper,
we are focused on topological spaces $\cK$ composed of finitely
many $0$-cells or vertices, together with finitely many 1-cells  --- copies of the open interval $(0,1)$
having  $0$-cells of $\cK$ as their boundaries --- and finitely many 2-cells --- open disks having boundary in
the space formed by the union of the $0$-cells and $1$-cells of $\cK$. In the case where there are no
$2$-cells and the space $\cK$ is connected, the fundamental group of $\cK$ with respect to any $0$-cell
$x$ as basepoint coincides with the group $\pi(\G,x)$ defined by Lewis where graph $\G$ is obtained by
replacing each $1$-cell $C$ of $\cK$ by the pairs $(y,z)$ and $(z,y)$ where $y$ and $z$ form the
boundary of $C$ in $\cK$.

\begin{defn}
\label{Def:cKj}
Let $\G=(X,R)$ be a graph and $j\ge 0$. Define $\cK_j = \cK(\G,j)$ to be the cell complex with
\begin{itemize}
\item the vertex set $X$ of $\G$ as its $0$-cells;
\item one $1$-cell for each edge of $\G$: if
$\{y,z\} \in E(\G)$, there is a $1$-cell $C_{yz}$ with $y$ and
$z$ as its boundary points;
\item one $2$-cell for each cycle in $\G$ of length $j$ or less: If the closed walk $p=(x_0,x_1,\ldots,x_\ell=x_0$
is a cycle in $\G$ with $\ell \le j$, then we have a disk $C_p$ with boundary consisting of the vertices
and edges of $p$.
\end{itemize}
\end{defn}

Note that, for $i<j$, $\cK_i$ is a topological subspace of $\cK_j$ and, for $i<3$, $\cK_i$ can be identified with $\G$ itself,
viewed as a topological space.

\begin{thm}
Let $\G=(X,R)$ be a graph, let $x\in X$, and let $j\ge 0$. Then the fundamental group $\pi\left(\cK_j,x\right)$
of the cell complex $\cK_j$ as defined in  Definition \ref{Def:cKj} is isomorphic to $\pi(\G,x)/\pi(\G,x,j)$.
\end{thm}

\begin{proof}
View $\Gamma$ as a $1$-dimensional cell complex with vertex set $X$ and one
$1$-cell for each edge of $\Gamma$. Then the fundamental group of this $1$-complex
at $x$ is $\pi(\Gamma,x)$.

By Definition 2.5, the complex $\cK_j$ is obtained from $\Gamma$ by attaching one $2$-cell along each cycle in $\Gamma$ of length at most $j$. Let $\mathcal C_j$ be the set of such cycles. For each $p\in\mathcal C_j$, let
$$
\alpha_p:S^1\to\cK_1
$$
be the loop in $\cK_1(=\G)$ determined by $p$. Thus $\cK_j$ is obtained from $\cK_1$ by attaching a copy $D_p^2$ of the $2$-disk along $\alpha_p$ for each $p\in\mathcal C_j$.

We first consider the effect of attaching a single $2$-cell. Let
$$
Y=\cK_1 \cup_\alpha D^2,
$$
where $\alpha:S^1\to\cK_1$ is a loop. Choose $V$ to be an open neighborhood of the attached disk $D^2$ together with a thin neighborhood of its boundary circle, and let $U$ be the union of $\cK_1$ with the same neighborhood. Then $U\cup V=Y$, $U$ deformation retracts onto $\Gamma$, $V$ is contractible, and $U\cap V$ deformation retracts onto $S^1$. Therefore
$$
\pi(U,x)\cong \pi(\Gamma,x), \qquad \pi(V,x)\cong \{e\}, \qquad \pi(U\cap V,x)\cong \pi(S^1)$$

Let
$$
i_1:U\cap V\to U,\qquad i_2:U\cap V\to V
$$
be the inclusion maps, and let
$$
i_{1*}:\pi(U\cap V,x)\to \pi(U,x),\qquad i_{2*}:\pi(U\cap V,x)\to \pi(V,x)
$$
be the induced homomorphisms. By van Kampen's Theorem,
$$
\pi(Y,x)\cong
\bigl(\pi(U,x)\ast \pi(V,x)\bigr)\Big/
\Big\langle
i_{1*}(\omega)i_{2*}(\omega)^{-1}
:
\omega\in\pi(U\cap V,x)
\Big\rangle.
$$
Since $\pi(V,x)$ is trivial, $i_{2*}(\omega)=e$ for every $\omega\in\pi(U\cap V,x)$, and so
$$
\pi(Y,x)\cong
\pi(U,x)\Big/
\Big\langle
i_{1*}(\omega)
:
\omega\in\pi(U\cap V,x)
\Big\rangle.
$$

Now $\pi(U\cap V,x)\cong \mathbb Z$, so it is generated by the class of a loop going once around the boundary circle. Under the map $i_{1*}$, this generator maps to the class $[\alpha]$ in $\pi(U,x)$. Thus
$$
\pi(Y,x)\cong \pi(U,x)/\langle[\alpha]\rangle.
$$
Since $\pi(U,x)\cong \pi(\Gamma,x)$, it follows that
$$
\pi(Y,x)\cong \pi(\Gamma,x)/\langle[\alpha]\rangle.
$$

Applying this argument successively to all $2$-cells attached in $K_j$, we obtain
$$
\pi(K_j,x)\cong \pi(\Gamma,x)\Big/\Big\langle[\alpha_p]:p\in\mathcal C_j\Big\rangle.
$$
Therefore,
$$
\pi(K_j,x)\cong \pi(\Gamma,x)/\pi(\Gamma,x,j),
$$
\end{proof}

\begin{cor}
\label{Cor:simply}
Let $\G=(X,R)$ be a graph, let $x\in X$, let $j\ge 0$ and let $\cK_j$ be the cell complex $\cK_j$ as defined in  Definition \ref{Def:cKj}.
\begin{itemize}
\item For $j<3$, $\pi\left(\cK_j,x\right) \cong \pi(\G,x)$;
\item $\cK_j$ is simply connected if and only if $\pi(\G,x,j)=\pi(\G,x)$.
\end{itemize}
\end{cor}

\begin{remark}
Our main objective, given $\G$,  is to find the smallest integer $j$ for which $\pi(\G,x,j)=\pi(\G,x)$.
In view of these results, the relationship between $\pi(\G,x,j)$ and $\pi(\G,x)$ is best described
in terms of the fundamental group of the cell complex $\cK_j$.
\end{remark}

\begin{remark}
One may use van Kampen's Theorem to verify the well-known fact that the fundamental group of any
simplicial complex is the same as that of its 2-skeleton. In our setting, this implies that, whenever
$\pi(\G,x,3)=\pi(\G,x)$, the \emph{clique complex} of $\G$, the simplicial complex formed by attaching
$(m-1)$-cells onto all $m$ cliques for all $m$, is simply connected.
\end{remark}

\section{Cost Functions}
\label{Subsec:Cost}

Before we address specific families of graphs, let us set up some general notation that will guide
our approach by mathematical induction.  Let $\G=(X,R)$ be a graph with path length distance
$\partial(x,y)$ on $X$.  Fix $x\in X$ and let $p\in \psi(\G,x)$. Write
$p = (x=x_0,x_1,\dots,x_l=x)$,  a closed walk  in $\G$ based at $x$ of length $l$. While the \emph{length}
of $p$ is the number of edges in the walk, and we may write $l(p)=l$, we also define the
\emph{cost} of $p$  (with respect to $x$) to be
\begin{equation}
\label{Eqn:cost}
c(p) = \sum_{ 1\leq i \leq l} \partial(x,x_i).
\end{equation}
Observe that $c(p)\ge 0$ and $c(p)=0$ if and only if $p = (x)$. For a fixed $k \ge 3$, our
general strategy is to find closed walks $q,p'\in \psi(\G,x)$ such that
\begin{itemize}
\item $p \sim q p'$
\item $q$ has essential length at most $k$
\item $c(p') < c(p)$ or $c(p) \le k$.
\end{itemize}
In Section \ref{Sec:Jnk}, we first describe this reduction procedure when $\G$ is a Johnson graph. The folded
Johnson graph is discussed next. In Section \ref{Sec:Hnq}, we consider the Hamming graphs and
folded cubes. We finish with some directions for further investigation.

\begin{example}
To illustrate this notation, we examine the trivial case where $\Gamma$ is a complete graph. After fixing a base vertex $x$, we have for all $y\in \Gamma$ with $y\ne x$ that $\partial(x,y)=1$, and so the cost of a walk $p$ is simply its length, $c(p)=l(p)$. Then the problem of decomposing $p$ into walks $q$ and $p'$ with $c(p')<c(p)$ requires that $l(p')<l(p)$.
\end{example}

\section{Johnson Graphs}
\label{Sec:Jnk}

We now introduce our notation for the Johnson graph $J(n,k)$ \cite[Section 9.1]{bcn}.
Let $N=[n]:=\{1,2,\dots, n\}$.  The vertex set $X$ of $J(n,k)$ is the set of all $k$-subsets
of $N$, denoted $N \choose k$.
For any $x,y\in X$, $\partial(x,y)=k-|x\cap y|$ so that $x$ and $y$ are
adjacent precisely when $\vert x\cap y \vert=k-1$. The goal of this section is to prove the following.

\begin{thm}
\label{Thm:Jnk}
Let $n$ and $k$ be positive integers with $n\ge 2k$, and let $\G = J(n,k)$. Then,
for any basepoint $x$, $\pi(\G,x)$ is generated by equivalence classes of walks of
essential length three; that is $\pi(\G,x,3)=\pi(\G,x)$.
Consequently,  the cell complex $\cK_3= \cK(\G,3)$ is simply connected.
\end{thm}

\begin{proof}
Let $\G = J(n,k)$, and let $x\in X$. Since the Johnson graph is vertex-transitive,
we may assume without loss of generality that $x = [k]=\{1,2,3,\dots,k\}$; this will be our choice
of basepoint throughout. For $p\in \psi(\G,x)$, let $\ell = \max \{\partial(x,x_i):x_i\in p\}$, the largest distance from $x$ achieved by any vertex in $p$); let
$m = \min \{i:\partial(x,x_i) = \ell \}$, the first index $m$ such that $x_m$ achieves
the maximum distance from $x$.

Clearly if $l(p)\le 3$ then $[p]\in \pi(\G,x,3)$, so we consider $p\in \psi(\G,x)$ with $l(p)\ge 4$. The goal of the induction process is to decompose $p$ into $p = q_1 q_2 \cdots q_r$ up to homotopy where $q_1,\dots,q_r$ each have essential length at most 3.

Using the cost function (\ref{Eqn:cost}), at each step of the induction we decompose a walk $p_i\in \psi(\G,x)$ into a product of walks $q_1\cdots q_kp_{i+1}$, each in $\psi(\G,x)$, such that $q_1,\dots,q_k$ each have essential length 3 and $c(p_{i+1})<c(p_i)$ when $i\ge 1$. Then we repeat the process with $p_{i+1}$. The process terminates at step $t$ when $c(p_t)\le 3$, as this implies $l(p_t) \le 3$.

In each case below, we have a closed walk
$p=(x=x_0,x_1,\ldots,x_{m-1},x_m,x_{m+1},\ldots,x_l = x)$ with subwalks
$r= (x_0,x_1,\ldots,x_{m-1})$ and $s=(x_{m+1},\ldots,x_l)$.

\bigskip

\noindent{\bf Case 0:} $x_{m-1}=x_{m+1}$ \\
Let $z=x_{m-1}=x_{m+1}$. In this case we have $p = r (z,x_m,z) s$.
Note that $p_0=r s$
is homotopic to $p$, so we may replace $p$ with $p_0$ with cost $c(p_0)=c(p)- \partial(x,x_m)-\partial(x,z)<c(p)$.

\bigskip

\noindent{\bf Case 1:} $x_{m-1}$ and $x_{m+1}$ are adjacent \\
One can write $p \sim q_0 p_0$ where $q_0 =  r (x_{m-1},x_m,x_{m+1},x_{m-1}) r^{-1}$ and
$p_0 = r ( x_{m-1},x_{m+1}) s$.
The walk $q_0$ has essential length 3: it cyclically reduces to the cycle $(x_{m-1}, x_m, x_{m+1}, x_{m-1})$. And we have $c(p_0) = c(p) -\ell < c(p)$.

\begin{figure} [h]
\resizebox{2.1in}{!}{
\begin{tikzpicture}[xscale=1,yscale=0.75,
    vert/.style={draw,fill=black,circle,inner sep=1.5pt}]
    \node[vert, blue] (x) at (-1,-0.5) {}; \node (lx) at (-1.4,-0.5) {$x$};
    \node[vert, blue] (xm-1) at (1,1) {}; \node (lxm-1) at (1,1.5) {$x_{m-1}$};
    \node[vert, blue] (xm) at (2,0) {}; \node (lxm) at (2.5,0) {$x_{m}$};
    \node[vert, blue] (xm+1) at (2,-1.5) {}; \node (lxm+1) at (2,-2) {$x_{m+1}$};

    \path (x) [thick, bend left=30, blue, dashed, -{Stealth[length=3mm]}] edge node [auto] {} (xm-1);
    \draw [thick, blue, -{Stealth[length=3mm]}]
    (xm-1) -- (xm);
    \draw [thick, blue, -{Stealth[length=3mm]}]
    (xm) -- (xm+1);
    \path (xm+1) [thick, bend left=30, blue, dashed, -{Stealth[length=3mm]}] edge node [auto] {} (x);

    \draw
    (xm-1) -- (xm+1);

\end{tikzpicture}
}
\resizebox{2.1in}{!}{
\begin{tikzpicture}[xscale=1,yscale=0.75,vert/.style={draw,fill=black,circle,inner sep=1.5pt}]
    \node[vert, red] (x) at (-1,-0.5) {}; \node (lx) at (-1.4,-0.5) {$x$};
    \node[vert, red] (xm-1) at (1,1) {}; \node (lxm-1) at (1,1.5) {$x_{m-1}$};
    \node[vert, red] (xm) at (2,0) {}; \node (lxm) at (2.5,0) {$x_{m}$};
    \node[vert, red] (xm+1) at (2,-1.5) {}; \node (lxm+1) at (2,-2) {$x_{m+1}$};

    \path (x) [thick, bend left=30, red, dashed, -{Stealth[length=3mm]}] edge node [auto] {} (xm-1);
    \draw [thick, red, -{Stealth[length=3mm]}]
    (xm-1) -- (xm);
    \draw [thick, red, -{Stealth[length=3mm]}]
    (xm) -- (xm+1);
    \draw [thick, red, -{Stealth[length=3mm]}]
    (xm+1) -- (xm-1);
    \path (xm-1) [thick, bend right=20, red, dashed, -{Stealth[length=3mm]}] edge node [auto] {} (x);

    \path (xm+1) [bend left=30, dashed] edge node [auto] {} (x);

\end{tikzpicture}

}
\resizebox{2.1in}{!}{
\begin{tikzpicture}[xscale=1,yscale=0.75,vert/.style={draw,fill=black,circle,inner sep=1.5pt},
black_vert/.style={draw,fill=black,circle,inner sep=1.25pt}]]

    \node[vert, teal] (x) at (-1,-0.5) {}; \node (lx) at (-1.4,-0.5) {$x$};
    \node[vert, teal] (xm-1) at (1,1) {}; \node (lxm-1) at (1,1.5) {$x_{m-1}$};
    \node[black_vert] (xm) at (2,0) {}; \node (lxm) at (2.5,0) {$x_{m}$};
    \node[vert, teal] (xm+1) at (2,-1.5) {}; \node (lxm+1) at (2,-2) {$x_{m+1}$};

    \path (x) [thick, bend left=30, teal, dashed, -{Stealth[length=3mm]}] edge node [auto] {} (xm-1);
    \draw [thick, teal, -{Stealth[length=3mm]}]
    (xm-1) -- (xm+1);
    \path (xm+1) [thick, bend left=30, teal, dashed, -{Stealth[length=3mm]}] edge node [auto] {} (x);

    \draw
    (xm-1) -- (xm)
    (xm) -- (xm+1);

\end{tikzpicture}
}
\caption{\label{Fig:3walks} {Case 1 in $J(n,k)$. The walk $p$ is in blue, $q_0$ in red, and $p_0$ in teal.}}
\end{figure}
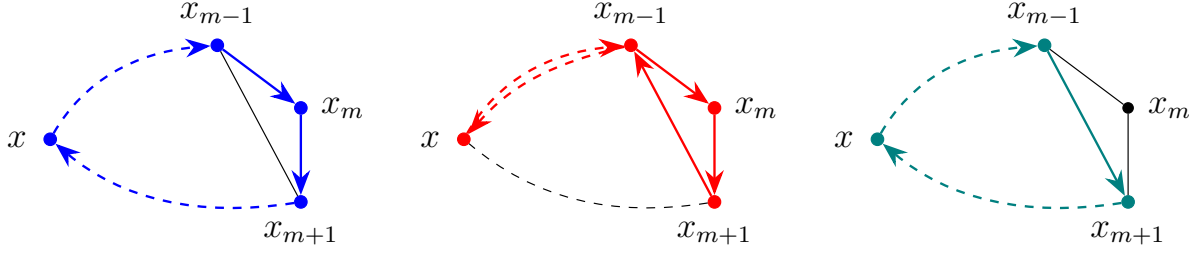
\bigskip

\noindent{\bf Case 2:} $x_{m-1}$ and $x_{m+1}$ are not adjacent (so $\partial(x_{m-1},x_{m+1})=2$) \\
If $\ell=1$ then $x_m=x_1$ and we have $x_{m-1}=x$ adjacent to $x_{m+1}=x_2$, which falls under Case 1. Therefore we may assume that $\ell \ge 2$. Notice that $\partial(x,x_{m-1})=\ell-1$, and so there is one element of the $(k-\ell+1)$-set $x\cap x_{m-1}$ which is not in $x_m$. Without loss of generality, assume $1\in x\cap x_{m-1}$ but $1\notin x_m$.  Since $\partial(x_{m-1},x_{m+1})=2$, the set
$S=x_{m-1}\cap x_{m+1}$ has  $\vert S\vert = k-2$. So we can find distinct elements $a,b,c\in N$ such that
\[
x_{m-1} = \{1,a\} \cup S, \qquad x_m = \{a,b\} \cup S, \qquad x_{m+1} = \{b,c\} \cup S.
\]
Note that the induced subgraph of $\G$ on the vertex set
$\{ z : \partial(z,x_{m-1})+ \partial(z,x_{m+1})=2 \}$ is isomorphic to the Johnson graph $J(4,2)$, the octahedron. Consider the two vertices in the set $\G(x_{m-1})\cap \G(x_{m})\cap \G(x_{m+1})$, say $y=\{1,b\}\cup S$ and $z = \{a, c\} \cup S$. Using $\hat{w}$ to denote the 2-set $w \backslash S$, these are depicted in Figure
\ref{Fig:octahedron}.

\begin{figure} [h]
\begin{center}
\begin{tikzpicture}[xscale=1,yscale=1,
    vert/.style={draw,fill=black,circle,inner sep=1.5pt}]

    \node[vert] (1a) at (0,2) {}; \node (l1a) at (0,2.5) {$\hat{x}_{m-1} = 1a$};
    \node[vert] (1c) at (-2,-0.25) {}; \node (l1c) at (-2.5,-0.25) {$1c$};
    \node[vert] (ac) at (-0.3,0.3) {};
    \node[vert] (1b) at (0.3,-0.3) {};
    \node[vert] (ab) at (2,0.25) {}; \node (lab) at (3,0.25) {$\hat{x}_m = ab$};
    \node[vert] (bc) at (0,-2) {}; \node (1bc) at (0,-2.5) {$\hat{x}_{m+1} = bc$};

    \draw[thick]
        (1a) -- (ab) -- (bc) -- (1b) -- (1a)
        (1c) -- (ac) -- (ab) -- (1b) -- (1c)
        (1a) -- (1c) -- (bc) -- (ac) -- (1a);

    \node (lac) at (-1.1,0.5) {$\hat{z} \ =ac$};
    \node (l1b) at (1.15,-0.6) {$\hat{y}= \ 1b$};

\end{tikzpicture}
\caption{\label{Fig:octahedron} For any two vertices $x,y$ with $\partial(x,y)=2$, the induced subgraph on $x$, $y$, and their common neighbors, is an octahedron. The vertex $\{k,\ell\}\cup S$ is abbreviated $k\ell$.}
\end{center}
\end{figure}
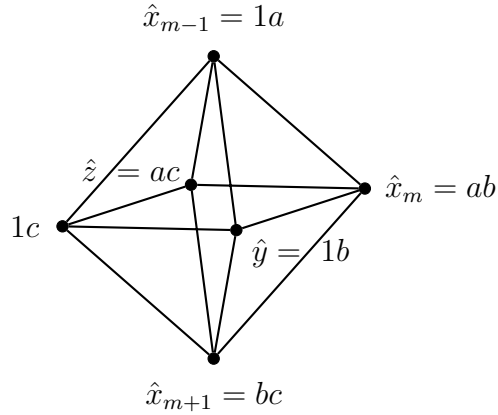

Using these common neighbors, we can write $p \sim q_0q_1p_0$ where
\begin{eqnarray*}
q_0 &=& r  (x_{m-1},x_m,y,x_{m-1}) r^{-1} \, , \\
q_1 &=& r (x_{m-1},y,x_m,x_{m+1},y,x_{m-1})  r^{-1} \, , \\
p_0 &=& r (x_{m-1},y,x_{m+1}) s ~ .
\end{eqnarray*}
Observe that $q_0$ and $q_1$ each have essential length 3, as $q_0 = r q_0'r^{-1}$  with
$q_0'=(x_{m-1},x_m,y,x_{m-1})$ and $q_1 = r_1q_1'r_1^{-1}$ with
$r_1 = (x,x_1,\dots, x_{m-1},y)$ and $q_1' = (y,x_m,x_{m+1},y)$.
The cost of $p_0$ is strictly less than the cost of $p$:
$$c(p) = \left( \sum_{\substack{1\leq i\leq n\\i\ne m}} \partial(x,x_i) \right) + \ell$$
$$c(p_0) = \left( \sum_{\substack{1 \leq i\leq n\\|i-m|>1}} \partial(x,x_i) \right) + \ell-1+\ell-1+\partial(x,x_{m+1}) = c(p)-1.$$

In all cases, we have shown that any closed walk $p$ based at $x$ with $c(p)>3$ can be expressed, up to homotopy, as a concatenation of walks all but one of which having essential length three and the last having
cost less than $c(p)$. This completes the proof, by induction.
\end{proof}

\bigskip
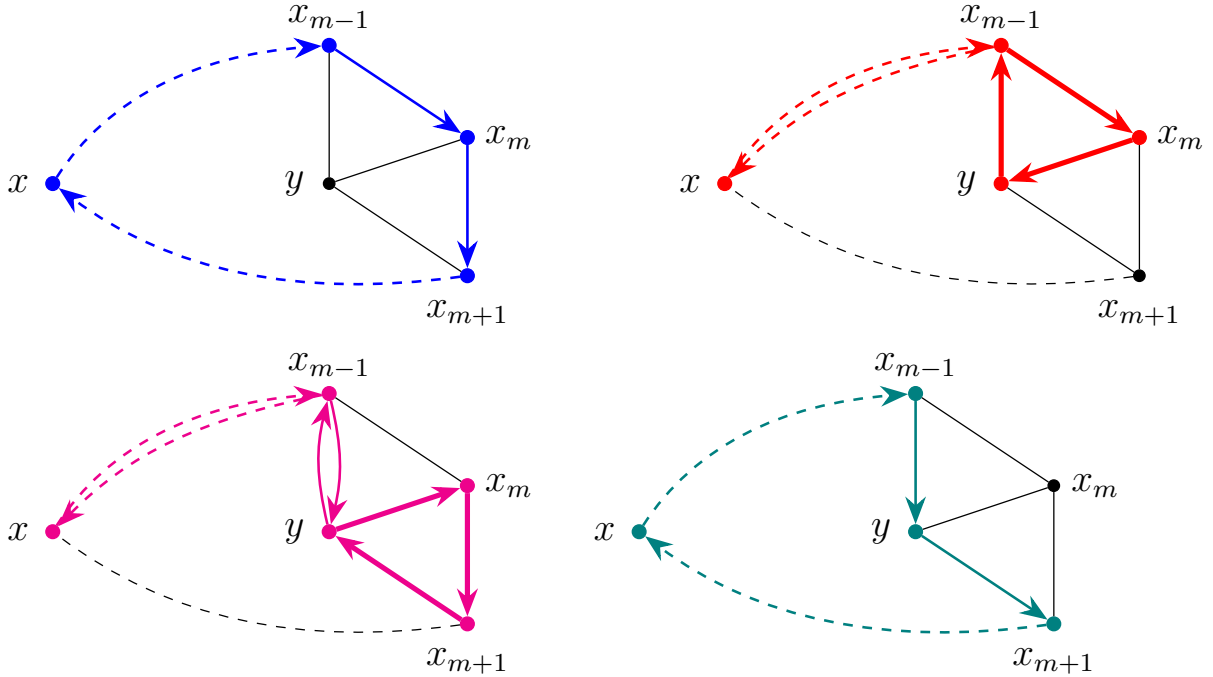
\begin{figure} [h]
\resizebox{3in}{!}{
\begin{tikzpicture}[xscale=1.5,yscale=1,
    blue_vert/.style={draw,blue, fill=blue,circle,inner sep=1.5pt},
    black_vert/.style={draw,fill=black,circle,inner sep=1.25pt}]

    \node[blue_vert] (x) at (-1,-0.5) {}; \node (lx) at (-1.25,-0.5) {$x$};
    \node[blue_vert] (xm-1) at (1,1) {}; \node (lxm-1) at (1,1.3) {$x_{m-1}$};
    \node[blue_vert] (xm) at (2,0) {}; \node (lxm) at (2.3,0) {$x_{m}$};
    \node[blue_vert] (xm+1) at (2,-1.5) {}; \node (lxm+1) at (2,-1.9) {$x_{m+1}$};
    \node[black_vert] (y) at (1,-0.5) {}; \node (ly) at (0.75,-0.5) {$y$};

    \draw
    (y) -- (xm-1)
    (y) -- (xm)
    (y) -- (xm+1);

    \path (x) [thick, bend left=30, blue, dashed, -{Stealth[length=3mm]}] edge node [auto] {} (xm-1);
    \draw [thick, blue, -{Stealth[length=3mm]}]
    (xm-1) -- (xm);
    \draw [thick, blue, -{Stealth[length=3mm]}]
    (xm) -- (xm+1);
    \path (xm+1) [thick, bend left=30, blue, dashed, -{Stealth[length=3mm]}] edge node [auto] {} (x);

\end{tikzpicture}
}
\resizebox{3in}{!}{
\begin{tikzpicture}[xscale=1.5,yscale=1,
    red_vert/.style={draw,red, fill=red,circle,inner sep=1.5pt},
    black_vert/.style={draw,fill=black,circle,inner sep=1.25pt}]

    \node[red_vert] (x) at (-1,-0.5) {}; \node (lx) at (-1.25,-0.5) {$x$};
    \node[red_vert] (xm-1) at (1,1) {}; \node (lxm-1) at (1,1.3) {$x_{m-1}$};
    \node[red_vert] (xm) at (2,0) {}; \node (lxm) at (2.3,0) {$x_{m}$};
    \node[black_vert] (xm+1) at (2,-1.5) {}; \node (lxm+1) at (2,-1.9) {$x_{m+1}$};
    \node[red_vert] (y) at (1,-0.5) {}; \node (ly) at (0.75,-0.5) {$y$};

    \draw
    (xm) -- (xm+1)
    (y) -- (xm+1);

    \path (xm+1) [bend left=30, dashed] edge node [auto] {} (x);

    \path (x) [thick, bend left=30, red, dashed, -{Stealth[length=3mm]}] edge node [auto] {} (xm-1);
    \path (xm-1) [thick, bend right=20, red, dashed, -{Stealth[length=3mm]}] edge node [auto] {} (x);
    \draw [ultra thick, red, -{Stealth[length=3mm]}]
    (xm-1) -- (xm);
    \draw [ultra thick, red, -{Stealth[length=3mm]}]
    (xm) -- (y);
    \draw [ultra thick, red, -{Stealth[length=3mm]}]
    (y) -- (xm-1);

\end{tikzpicture}
}
\resizebox{3in}{!}{
\begin{tikzpicture}[xscale=1.5,yscale=1,
    mag_vert/.style={draw,magenta, fill=magenta,circle,inner sep=1.5pt},
    black_vert/.style={draw,fill=black,circle,inner sep=1.25pt}]

    \node[mag_vert] (x) at (-1,-0.5) {}; \node (lx) at (-1.25,-0.5) {$x$};
    \node[mag_vert] (xm-1) at (1,1) {}; \node (lxm-1) at (1,1.3) {$x_{m-1}$};
    \node[mag_vert] (xm) at (2,0) {}; \node (lxm) at (2.3,0) {$x_{m}$};
    \node[mag_vert] (xm+1) at (2,-1.5) {}; \node (lxm+1) at (2,-1.9) {$x_{m+1}$};
    \node[mag_vert] (y) at (1,-0.5) {}; \node (ly) at (0.75,-0.5) {$y$};

    \draw
    (xm-1) -- (xm);

    \path (xm+1) [bend left=30, dashed] edge node [auto] {} (x);

    \path (x) [thick, bend left=30, magenta, dashed, -{Stealth[length=3mm]}] edge node [auto] {} (xm-1);
    \path (xm-1) [thick, bend right=20, magenta, dashed, -{Stealth[length=3mm]}] edge node [auto] {} (x);

    \path (xm-1) [thick, bend left=9, magenta, -{Stealth[length=3mm]}] edge node [auto] {} (y);
    \path (y) [thick, bend left=9, magenta, -{Stealth[length=3mm]}] edge node [auto] {} (xm-1);

    \draw [ultra thick, magenta, -{Stealth[length=3mm]}]
    (y) -- (xm);
    \draw [ultra thick, magenta, -{Stealth[length=3mm]}]
    (xm) -- (xm+1);
    \draw [ultra thick, magenta, -{Stealth[length=3mm]}]
    (xm+1) -- (y);

\end{tikzpicture}
}
\resizebox{3in}{!}{
\begin{tikzpicture}[xscale=1.5,yscale=1,
    teal_vert/.style={draw,teal, fill=teal,circle,inner sep=1.5pt},
    black_vert/.style={draw,fill=black,circle,inner sep=1.25pt}]

    \node[teal_vert] (x) at (-1,-0.5) {}; \node (lx) at (-1.25,-0.5) {$x$};
    \node[teal_vert] (xm-1) at (1,1) {}; \node (lxm-1) at (1,1.3) {$x_{m-1}$};
    \node[black_vert] (xm) at (2,0) {}; \node (lxm) at (2.3,0) {$x_{m}$};
    \node[teal_vert] (xm+1) at (2,-1.5) {}; \node (lxm+1) at (2,-1.9) {$x_{m+1}$};
    \node[teal_vert] (y) at (1,-0.5) {}; \node (ly) at (0.75,-0.5) {$y$};

    \draw
    (xm-1) -- (xm)
    (y) -- (xm)
    (xm) -- (xm+1);

    \path (x) [thick, bend left=30, teal, dashed, -{Stealth[length=3mm]}] edge node [auto] {} (xm-1);
    \draw [thick, teal, -{Stealth[length=3mm]}]
    (xm-1) -- (y);
    \draw [thick, teal, -{Stealth[length=3mm]}]
    (y) -- (xm+1);
    \path (xm+1) [thick, bend left=30, teal, dashed, -{Stealth[length=3mm]}] edge node [auto] {} (x);

\end{tikzpicture}
}
\caption{\label{Fig:Johnsoncase2}Case 2 in $J(n,k)$. Walk $p$ is in blue, $q_0$ in red, $q_1$ in magenta, and $p_0$ in teal.}
\end{figure}

\begin{example} \label{Ex:J63} We apply the induction process to a walk of length 7 in $J(6,3)$ with basepoint $x=\{1,2,3\}$, depicted in Figure \ref{Fig:J63}. The vertex $\{a,b,c\}$ is denoted $abc$.

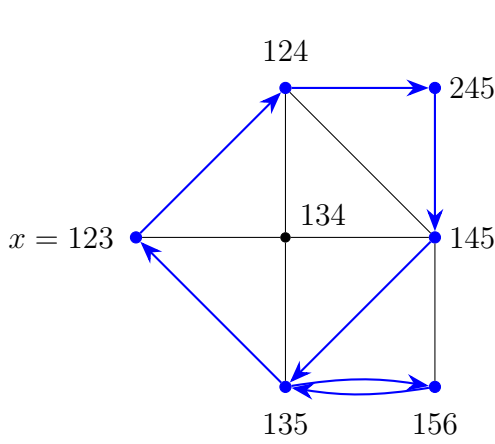
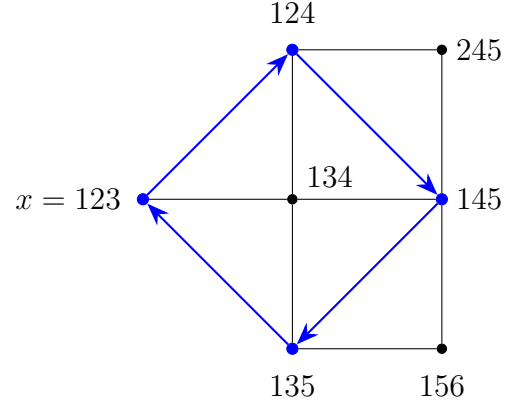
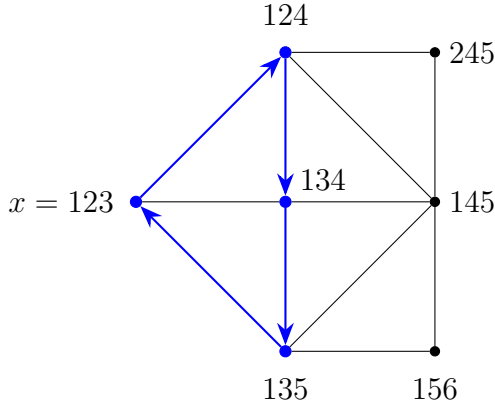
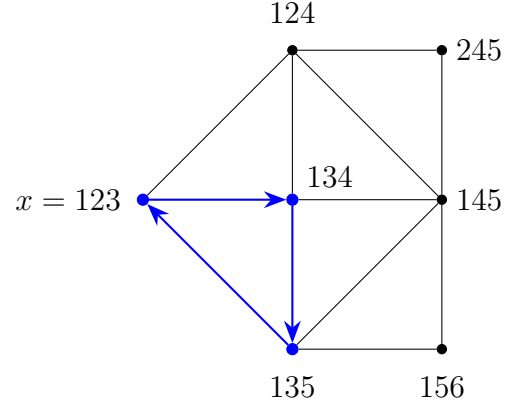
\begin{figure}[h]
\begin{subfigure}[b]{0.5\textwidth}
\resizebox{2.75in}{!}{
\begin{tikzpicture}[xscale=1,yscale=1, blue_vert/.style={draw,blue,fill=blue,circle,inner sep=1.5pt}, vert/.style={draw,fill=black,circle,inner      sep=1.25pt}]

    \node[blue_vert] (123) at (0,0) {}; \node (lx) at (-1,0,0) {$x=123$};
    \node[blue_vert] (124) at (2,2) {}; \node (lx1) at (2,2.5) {$124$};
    \node[blue_vert] (245) at (4,2) {}; \node (lx2) at (4.5,2) {$245$};
    \node[blue_vert] (145) at (4,0) {}; \node (lx2) at (4.5,0) {$145$};
    \node[vert] (134) at (2,0) {}; \node (lx2) at (2.5,0.3) {$134$};
    \node[blue_vert] (135) at (2,-2) {}; \node (lx2) at (2,-2.5) {$135$};
    \node[blue_vert] (156) at (4,-2) {}; \node (lx2) at (4,-2.5) {$156$};

    \draw[thick, blue, -{Stealth[length=3mm]}]      (123) -- (124);
    \draw[thick, blue, -{Stealth[length=3mm]}]      (124) -- (245);
    \draw[thick, blue, -{Stealth[length=3mm]}]      (245) -- (145);
    \draw[thick, blue, -{Stealth[length=3mm]}]      (145) -- (135);
    \path (135) [thick, bend left=9, blue, -{Stealth[length=3mm]}] edge node [auto] {} (156);
    \path (156) [thick, bend left=9, blue, -{Stealth[length=3mm]}] edge node [auto] {} (135);
    \draw[thick, blue, -{Stealth[length=3mm]}]      (135) -- (123);

    \draw
    (123) -- (134)
    (124) -- (134)
    (134) -- (135)
    (124) -- (145)
    (134) -- (145)
    (145) -- (156);

\end{tikzpicture}
}
\caption{The original walk.}
\end{subfigure}
\begin{subfigure}[b]{0.5\textwidth}
\resizebox{2.75in}{!}{
\begin{tikzpicture}[xscale=1,yscale=1, blue_vert/.style={draw,blue,fill=blue,circle,inner sep=1.5pt}, vert/.style={draw,fill=black,circle,inner      sep=1.25pt}]

    \node[blue_vert] (123) at (0,0) {}; \node (lx) at (-1,0,0) {$x=123$};
    \node[blue_vert] (124) at (2,2) {}; \node (lx1) at (2,2.5) {$124$};
    \node[vert] (245) at (4,2) {}; \node (lx2) at (4.5,2) {$245$};
    \node[blue_vert] (145) at (4,0) {}; \node (lx2) at (4.5,0) {$145$};
    \node[vert] (134) at (2,0) {}; \node (lx2) at (2.5,0.3) {$134$};
    \node[blue_vert] (135) at (2,-2) {}; \node (lx2) at (2,-2.5) {$135$};
    \node[vert] (156) at (4,-2) {}; \node (lx2) at (4,-2.5) {$156$};

    \draw[thick, blue, -{Stealth[length=3mm]}]      (123) -- (124);
    \draw[thick, blue, -{Stealth[length=3mm]}]      (124) -- (145);
    \draw[thick, blue, -{Stealth[length=3mm]}]      (145) -- (135);
    \draw[thick, blue, -{Stealth[length=3mm]}]      (135) -- (123);

    \draw
    (123) -- (134)
    (124) -- (134)
    (134) -- (135)
    (124) -- (245)
    (134) -- (145)
    (135) -- (156)
    (245) -- (145)
    (145) -- (156);

\end{tikzpicture}
}
\caption{Case 0 removes 156 from the walk; case 1 removes 245.}
\end{subfigure}
\begin{subfigure}[b]{0.5\textwidth}
\resizebox{2.75in}{!}{
\begin{tikzpicture}[xscale=1,yscale=1, blue_vert/.style={draw,blue,fill=blue,circle,inner sep=1.5pt}, vert/.style={draw,fill=black,circle,inner      sep=1.25pt}]

    \node[blue_vert] (123) at (0,0) {}; \node (lx) at (-1,0,0) {$x=123$};
    \node[blue_vert] (124) at (2,2) {}; \node (lx1) at (2,2.5) {$124$};
    \node[vert] (245) at (4,2) {}; \node (lx2) at (4.5,2) {$245$};
    \node[vert] (145) at (4,0) {}; \node (lx2) at (4.5,0) {$145$};
    \node[blue_vert] (134) at (2,0) {}; \node (lx2) at (2.5,0.3) {$134$};
    \node[blue_vert] (135) at (2,-2) {}; \node (lx2) at (2,-2.5) {$135$};
    \node[vert] (156) at (4,-2) {}; \node (lx2) at (4,-2.5) {$156$};

    \draw[thick, blue, -{Stealth[length=3mm]}]      (123) -- (124);
    \draw[thick, blue, -{Stealth[length=3mm]}]      (124) -- (134);
    \draw[thick, blue, -{Stealth[length=3mm]}]      (134) -- (135);
    \draw[thick, blue, -{Stealth[length=3mm]}]      (135) -- (123);

    \draw
    (123) -- (134)
    (124) -- (245)
    (124) -- (145)
    (134) -- (145)
    (135) -- (145)
    (135) -- (156)
    (245) -- (145)
    (145) -- (156);

\end{tikzpicture}
}
\caption{Case 2 exchanges 145 for 134.}
\end{subfigure}
\begin{subfigure}[b]{0.5\textwidth}
\resizebox{2.75in}{!}{
\begin{tikzpicture}[xscale=1,yscale=1, blue_vert/.style={draw,blue,fill=blue,circle,inner sep=1.5pt}, vert/.style={draw,fill=black,circle,inner      sep=1.25pt}]

    \node[blue_vert] (123) at (0,0) {}; \node (lx) at (-1,0,0) {$x=123$};
    \node[vert] (124) at (2,2) {}; \node (lx1) at (2,2.5) {$124$};
    \node[vert] (245) at (4,2) {}; \node (lx2) at (4.5,2) {$245$};
    \node[vert] (145) at (4,0) {}; \node (lx2) at (4.5,0) {$145$};
    \node[blue_vert] (134) at (2,0) {}; \node (lx2) at (2.5,0.3) {$134$};
    \node[blue_vert] (135) at (2,-2) {}; \node (lx2) at (2,-2.5) {$135$};
    \node[vert] (156) at (4,-2) {}; \node (lx2) at (4,-2.5) {$156$};

    \draw[thick, blue, -{Stealth[length=3mm]}]      (123) -- (134);
    \draw[thick, blue, -{Stealth[length=3mm]}]      (134) -- (135);
    \draw[thick, blue, -{Stealth[length=3mm]}]      (135) -- (123);

    \draw
    (123) -- (124)
    (124) -- (134)
    (124) -- (245)
    (124) -- (145)
    (134) -- (145)
    (135) -- (145)
    (135) -- (156)
    (245) -- (145)
    (145) -- (156);

\end{tikzpicture}
}
\caption{Case 1 removes 124. The final walk is a triangle.}
\end{subfigure}
\caption{\label{Fig:J63} Reduction of a walk in $J(6,3)$.}
\end{figure}

\end{example}

\subsection{Folded Johnson graphs}
\label{Subsec:fJ2kk}

When $n=2k$, the Johnson graph $J(n,k)$ is antipodal: each vertex $x \in N\choose k$ is at distance $k$ from a unique vertex in $J(n,k)$, its antipode $N\backslash x$. Identifying antipodal pairs gives us another distance-regular graph, the \emph{folded Johnson graph}, $\bar{J}(2k,k)$ \cite[Section 9.1C]{bcn}. It is well-known that $\bar{J}(2k,k)$ has diameter $d = \lfloor k/2 \rfloor$ and is isomorphic to the graph $\G$ obtained by merging distances $1$ and $k-1$ in the Johnson graph $J(2k-1,k-1)$. Using this isomorphism, let us make the covering map explicit.  Let $Y = \binom{[2k]}{k}$
and let $X$ denote the vertex set of $J(2k-1,k-1)$. We have
\begin{equation}
\label{Eqn:JohnsonCover}
\kappa : Y \rightarrow X \qquad \text{via} \qquad
\kappa(x) = \begin{cases}
x \setminus \{n\}&\text{if $n\in x$\ ;} \\
[2k-1] \setminus x&\text{if $n\not\in x$\ .}
\end{cases}
\end{equation}
For $k\ge 4$, this is a local isomorphism where two elements $x,y \in X$ are adjacent in $\G$ if $|x\cap y|=k-2$ or $|x\cap y|= 0$. Let us call edges of the first type \emph{short} edges and edges of the second type \emph{long} edges. Note that $\kappa$ extends naturally to a function from the set of walks in $J(2k,k)$ to some subcollection of the walks in $\G$.

Let $\cK_3$ denote the cell complex with 0-cells and 1-cells given by vertices and edges of $\G$ and $2$-cells given by the triangles in $\G$.

\begin{lem}
\label{Lem:triangles}
Assume $k>3$. Each triangle in $\G$ contains 0 or 2 long edges. For any (closed) disk $B$ in
$\cK_3$, the boundary  of  $B$ contains an even number of long edges.
\end{lem}

\begin{proof}
Since $n < 3k$, no triangle has three long edges. Since $k-1>2$, no triangle has just one long edge.
Now consider a disk composed of triangles and inductively verify that the mod two sum of the number
of long edges on the boundary remains zero modulo 2 as the number of triangles increases.
\end{proof}

\begin{remark}
For $k=3$, $\G$ is a complete graph and there exist triangles having one long edge. In this case only, we have
$\pi(\G,x,3)=\pi(\G,x)$.
\end{remark}

Results of Lewis \cite[Theorem 26]{lewis} tell us that the cell complex $\cK_k$ associated with $\bar{J}(2k,k)$ is simply connected when $k$ is odd. We now extend this result.

\begin{thm}
\label{Thm:foldJ2kk}
Let $k>3$,  let $\G=(X,R)$ be the folded Johnson graph $\bar{J}(2k,k)$ and let $x\in X$.
Consider the corresponding cell complex $\cK_j$, $3\le j\le k$, from Definition \ref{Def:cKj}.
Then $\pi(\G,x,k)=\pi(\G,x)$. So, $\cK_k$ is simply connected while
$\pi\left(\cK_j,x\right) \cong \ints/2\ints$ for $3\le j < k$.
\end{thm}

\begin{proof}
Let $x_0=\{1,\ldots,k\} \in Y$ and let $x=\kappa(x_0)$, employing the isomorphism above and our map (\ref{Eqn:JohnsonCover}).
We define a closed walk $\bar{p}$ based at $x$ as follows.
For $1\le j \le k$ let $x_j =  \{ j+1,\ldots, j+k\}$ in
$J(2k,k)$ so that $x_k$ is the antipode of $x_0$. The image $\bar{p}$ of the walk $\bar{p}'=(x_0,x_1,\ldots,x_k)$ in
$\G'=J(2k,k)$ under the projection $\kappa$ is a closed walk (a cycle, in fact) in $\G$ based
at $x$ which contains just one long edge. Since any disk composed of $2$-cells in the cell
complex $\cK_3$ coming from $\G$ contains an even number of long edges along its boundary, this cycle  $\bar{p}$ is not contractible. So $[\bar{p}] \neq e$ in $\pi(\cK_3,x)$.

Next, suppose $3<j<k$.
Theorem \ref{Thm:Jnk} tells us that any $2$-cell in $\cK(\G',j)$ (cf. Definition \ref{Def:cKj}) has boundary
equal to the boundary of some disk composed entirely of triangles. So any homotopy
$p\sim \cdots \sim (x)$ in $\cK(\G',j)$ can be converted into a homotopy $p\sim \cdots \sim (e)$
in $\cK(\G',3)$. From the previous paragraph, it follows that  $[\bar{p}] \neq e$ in
$\pi(\cK_j,x)$ for $j<k$.

On the other hand, if $p$ and $r$ are any closed walks in $\G$ based at $x$ such that $[p]\neq e$ and $[r]\neq e$, we will now show that $[pr] = e$ in $\pi(\cK_3,x)$. It  then follows that $\pi(\cK_3,x) \cong \ints/2\ints$ whenever $k\ge 4$. That is, for $k\ge 4$, $\pi(\G,x,3)$ is a subgroup  of index two in $\pi(\G,x)$.

Let $pr$ be the concatenation of these two closed walks in $\G$ based at $x$. Decompose
$pr$ as the concatenation of walks
\begin{equation}
\label{Eqn:pqqq}
pr = q_1 q_2 \cdots q_s
\end{equation}
where each $q_j$ is a closed walk based at $x$  and does not contain $x$ as an internal vertex.
If $[q_j] \neq e$, we may lift $q_j$ to a walk in $J(2k,k)$ in two ways: either from $x_0$ to its antipode $x_k$ or from $x_k$ to $x_0$; otherwise $q_j$ is contractible in $\cK_3$.

Now suppose, for $q_i,q_j\in \psi(\G,x)$,
that $q_i$ and $q_j$ are not contractible in $\cK_3$, and neither contains $x$ as an internal vertex. Then
the concatenation $q_i q_j$ can be viewed as a projection $q_i q_j = \kappa(q'_i q'_j)$ of a closed walk in $J(2k,k)$
based at $x_0$ which begins with a walk $q'_i$ from $x_0$ to its antipode satisfying $\kappa(q'_i)=q_i$ followed
by a walk $q'_j$ from the antipode $x_k$ to $x_0$ satisfying $\kappa(q'_j)=q_j$.
Now we use Theorem \ref{Thm:Jnk} to conclude that $q'_i q'_j$ is null homotopic in the cell complex
$\cK\left( \Gamma' , 3\right)$. Since $k>3$, each $2$-cell in $\cK\left( \Gamma' , 3\right)$ projects,
via $\kappa$, to a $2$-cell in $\cK_3$, we
have that $[q_i q_j]=e$ in $\pi( \cK_3,x)$.  It follows that $[pr]=[\bar{p}]$ or $[pr]=e$.

Finally, in the cell complex $\cK_k$, each $q_j$ in (\ref{Eqn:pqqq})
is null homotopic giving us that $\pi( \cK_k,x)$ is the trivial group.
\end{proof}

\section{Hamming graphs}
\label{Sec:Hnq}
We define the Hamming graph $H(d,q)$ as follows. Let $N=\{0,1,\dots,q-1\}$. The vertex set $X$ of
$H(d,q)$ is $N^d$. Thus a vertex is an ordered $d$-tuple
$
u=u_1 \cdots u_d
$
with each $u_i\in N$. Two vertices are adjacent when they differ in exactly one coordinate.
The distance $\partial(u,v)$ is the number of coordinates in which $u$ and $v$ differ.

The goal of this section is to prove the following.

\begin{thm}
\label{Thm:Hdq}
Let $d$ and $q$ be positive integers with $q\ge 2$, and let $\G=H(d,q)$. Then, for any basepoint $x$,
$\pi(\G,x)$ is generated by equivalence classes of walks of essential length at most four; that is,
$
\pi(\G,x,4)=\pi(\G,x).
$
Consequently, the cell complex $\cK_4=\cK(\G,4)$ is simply connected.
\end{thm}

\begin{proof}
Since the Hamming graph is vertex-transitive, we may assume without loss of generality that
$
x=00\cdots 0 ~.
$
Let
$
p=(x=x_0,x_1,\ldots,x_l=x)
$
be a closed walk based at $x$. We argue by induction on the cost $c(p)$ defined in (\ref{Eqn:cost}).
 If $l(p) \le 4$, then $[p]\in \pi(\G,x,4)$.

Assume now that $l(p)>4$ and $c(p)>0$. Let
$
\ell=\max\{\partial(x,x_i):x_i\text{ in }p\}
$
and let
$
m=\min\{i:\partial(x,x_i)=\ell\}.
$
Thus $x_m$ is the first vertex of $p$ at maximum distance from the basepoint. Set
$
a=x_{m-1}, b=x_m, c=x_{m+1}.
$
Since $b$ is the first vertex at distance $\ell$ from $x$ and $a$ is adjacent to $b$, we have
$
\partial(x,a)=\ell-1.
$
By permuting coordinates and relabeling the nonzero symbols in each coordinate, we may assume that
$
b = 1 \cdots 110 \cdots 0
$
has exactly $\ell$ nonzero coordinates and that
$
a = 1 \cdots 100 \cdots 0
$
is obtained from $b$ by changing $b_\ell$ from $1$ to $0$.

\bigskip

\noindent{\bf Case 0:} $a=c$. \\
In this case the walk $p$ contains a subwalk of the form
$
(a,b,a).
$
Deleting this subwalk gives a closed walk
$
p_0=(x=x_0,x_1,\ldots,x_{m-1},x_{m+2},\ldots,x_l=x).
$
By Lewis's homotopy relation, $p\sim p'$, since $p$ is obtained from $p_0$ by an extension, and we have
$
c(p_0)<c(p).
$

\bigskip

\noindent{\bf Case 1:} $a$ and $c$ are adjacent. \\
Then
$
(a,b,c,a)
$
is a triangle in $\G$. Let
$
r=(x=x_0,x_1,\ldots,x_{m-1}=a)
$
and let
$
s=(c=x_{m+1},x_{m+2},\ldots,x_l=x).
$
Define $q_0=r(a,b,c,a)r^{-1}$ and $p_0=r(a,c)s$.
The walk $q_0$ clearly has essential length three, since
$(a,b,c,a)$ is a triangle in $\G$.
Moreover, $p\sim q_0p_0$. Indeed, after deleting subwalks of the form $(u,v,u)$ from $q_0p_0$, one obtains
$p$. We have
$
c(p_0)=c(p)-\partial(x,b)<c(p).
$

\bigskip
\noindent{\bf Case 2:} $a$ and $c$ are not adjacent. \\
In this case $\ell\ge 2$.
Since $\G$ is distance-transitive, after permuting coordinates and relabeling symbols, we may write
\begin{eqnarray*}
b &=& \underbrace{1\cdots 1}_{\ell-2} \,110 \cdots 0 , \\
a &=& \underbrace{1\cdots 1}_{\ell-2} \,100 \cdots 0 , \\
c &=& \underbrace{1\cdots 1}_{\ell-2} \,\beta 10 \cdots 0 ,
\end{eqnarray*}
where $\beta\ne 1$. Choose
\[
y = \underbrace{1\cdots 1}_{\ell-2} \,\beta 00 \cdots 0 ~ .
\]
Then $y$ is adjacent to both $a$ and $c$, and
$
\partial(x,y)\le \ell-1.
$

Again let
$
r=(x=x_0,x_1,\ldots,x_{m-1}=a)
$
and
$
s=(c=x_{m+1},x_{m+2},\ldots,x_l=x).
$
Define
$$
q_0=r(a,b,c,y,a)r^{-1}
$$
and
$$
p_0=r(a,y,c)s ~ .
$$
The walk $q_0$ has essential length four  since
$(a,b,c,y,a)$ is a $4$-cycle in $\G$.
Moreover, $p\sim q_0p_0$, again by deleting subwalks of the form $(u,v,u)$ from $q_0p_0$. Since
$
\partial(x,b)=\ell
$
and
$
\partial(x,y)\le \ell-1,
$
we have
$
c(p_0)\le c(p)-1<c(p).
$

In all cases, either $p$ is homotopic to a closed walk of smaller cost, or $p$ is homotopic to a product
$q_0p_0$, where $q_0$ has essential length at most four and $c(p_0)<c(p)$. Repeating this process,
we express $p$ as a product of closed walks of essential length at most four. Hence
$
\pi(\G,x,4)=\pi(\G,x).
$
By Corollary \ref{Cor:simply}, the cell complex $\cK_4$ is simply connected.
\end{proof}

\begin{figure} [h]
\resizebox{2.4in}{!}{
\begin{tikzpicture}[xscale=1,yscale=0.75,blue_vert/.style={draw,blue,fill=blue,circle,inner sep=1.5pt}, vert/.style={draw,fill=black,circle,inner      sep=1.25pt}]
    \node[blue_vert] (x) at (-1,-0.5) {}; \node (lx) at (-1.4,-0.5) {$x$};
    \node[blue_vert] (xm-1) at (1,1) {}; \node (lxm-1) at (1,1.5) {$a=x_{m-1}$};
    \node[blue_vert] (xm) at (2,0) {}; \node (lxm) at (2.8,0) {$b=x_{m}$};
    \node[blue_vert] (xm+1) at (2,-1.5) {}; \node (lxm+1) at (2,-2) {$c=x_{m+1}$};
    \node[vert] (y) at (1,-0.5) {}; \node (ly) at (0.6,-0.5) {$y$};

    \path (x) [thick, bend left=30, blue, dashed, -{Stealth[length=3mm]}] edge node [auto] {} (xm-1);
    \draw [thick, blue, -{Stealth[length=3mm]}]
    (xm-1) -- (xm);
    \draw [thick, blue, -{Stealth[length=3mm]}]
    (xm) -- (xm+1);
    \path (xm+1) [thick, bend left=30, blue, dashed, -{Stealth[length=3mm]}] edge node [auto] {} (x);

    \draw
    (xm-1) -- (y)
    (y) -- (xm+1);

\end{tikzpicture}
}
\resizebox{2in}{!}{
\begin{tikzpicture}[xscale=1,yscale=0.75,red_vert/.style={draw,red,fill=red,circle,inner sep=1.5pt}, vert/.style={draw,fill=black,circle,inner      sep=1.25pt}]
    \node[red_vert] (x) at (-1,-0.5) {}; \node (lx) at (-1.4,-0.5) {$x$};
    \node[red_vert] (xm-1) at (1,1) {}; \node (lxm-1) at (1,1.5) {$a$};
    \node[red_vert] (xm) at (2,0) {}; \node (lxm) at (2.3,0) {$b$};
    \node[red_vert] (xm+1) at (2,-1.5) {}; \node (lxm+1) at (2,-2) {$c$};
    \node[red_vert] (y) at (1,-0.5) {}; \node (ly) at (0.6,-0.5) {$y$};

    \path (x) [thick, bend left=30, red, dashed, -{Stealth[length=3mm]}] edge node [auto] {} (xm-1);
    \draw [thick, red, -{Stealth[length=3mm]}]
    (xm-1) -- (xm);
    \draw [thick, red, -{Stealth[length=3mm]}]
    (xm) -- (xm+1);
    \draw [thick, red, -{Stealth[length=3mm]}]
    (xm+1) -- (y);
    \draw [thick, red, -{Stealth[length=3mm]}]
    (y) -- (xm-1);
    \path (xm-1) [thick, bend right=20, red, dashed, -{Stealth[length=3mm]}] edge node [auto] {} (x);

    \path (xm+1) [bend left=30, dashed] edge node [auto] {} (x);

\end{tikzpicture}

}
\resizebox{2in}{!}{
\begin{tikzpicture}[xscale=1,yscale=0.75,teal_vert/.style={draw,teal,fill=teal,circle,inner sep=1.5pt}, vert/.style={draw,fill=black,circle,inner      sep=1.25pt}]

    \node[teal_vert] (x) at (-1,-0.5) {}; \node (lx) at (-1.4,-0.5) {$x$};
    \node[teal_vert] (xm-1) at (1,1) {}; \node (lxm-1) at (1,1.5) {$a$};
    \node[vert] (xm) at (2,0) {}; \node (lxm) at (2.3,0) {$b$};
    \node[teal_vert] (xm+1) at (2,-1.5) {}; \node (lxm+1) at (2,-2) {$c$};
    \node[teal_vert] (y) at (1,-0.5) {}; \node (ly) at (0.6,-0.5) {$y$};

    \path (x) [thick, bend left=30, teal, dashed, -{Stealth[length=3mm]}] edge node [auto] {} (xm-1);
    \draw [thick, teal, -{Stealth[length=3mm]}]
    (xm-1) -- (y);
    \draw [thick, teal, -{Stealth[length=3mm]}]
    (y) -- (xm+1);
    \path (xm+1) [thick, bend left=30, teal, dashed, -{Stealth[length=3mm]}] edge node [auto] {} (x);

    \draw
    (xm-1) -- (xm)
    (xm) -- (xm+1);

\end{tikzpicture}
}
\caption{\label{Fig:hammingcase2} {Case 2 in $H(d,q)$. The walk $p$ is in blue, $q_0$ in red, and $p_0$ in teal.}}
\end{figure}
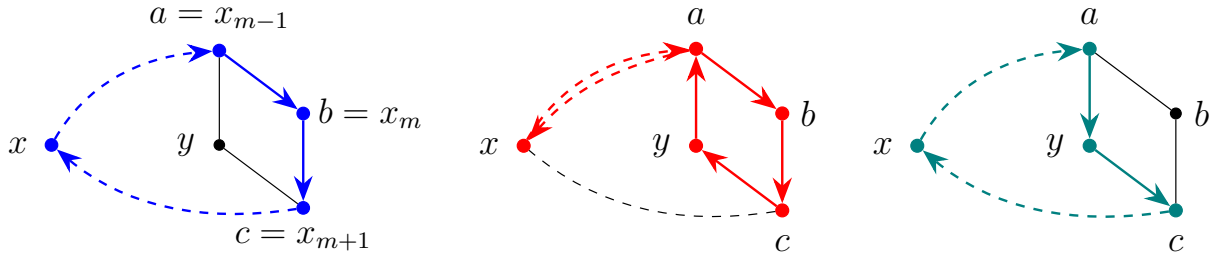
\bigskip

\subsection{Folded Cubes}
\label{Subsec:foldcube}

When $q=2$, the Hamming graph $H(n,q)$, also known as the $n$-cube and written $Q_n$,
is  both bipartite and antipodal. As with the Johnson graph $J(2k,k)$ above,  each vertex of
$H(n,2)$ is at distance $n$ from its unique antipode. Identifying antipodal pairs gives us the \emph{folded cube} \cite[Section 9.2D]{bcn}. It is well-known that this graph $\bar{H}(n,2)$ is isomorphic to the
graph obtained by merging distances $1$ and $n-1$ in the $(n-1)$-cube $H(n-1,2)$.  As in Section \ref{Subsec:fJ2kk}, we employ this isomorphism, denoting a bit flip by $b'=1-b$:
\begin{equation}
\label{Eqn:CubeCover}
\kappa : \{0,1\}^n \rightarrow \{0,1\}^{n-1}
\qquad \text{via} \qquad
\kappa(u_1 \cdots u_n) = \begin{cases}
u_2 \cdots  u_n &\text{if $u_1=0$;} \\
u_2' \cdots  u_n' &\text{if $u_1=1$.}
\end{cases}
\end{equation}
For $n>3$, this is a local isomorphism from the $n$-cube to the graph $\G$ formed by merging
distances $1$ and $n-1$ in the $(n-1)$-cube. For adjacent vertices $x$ and $y$ in $\G$, let us
say the edge $\{x,y\}$ is a \emph{short} edge
of $\G$ if $\partial(x,y)=1$ in the $(n-1)$-cube and that  $\{x,y\}$ is a \emph{long} edge
of $\G$ if $\partial(x,y)=n-1$ in the $(n-1)$-cube.  Just as before, $\kappa$ extends naturally
to a function from the set of walks in $H(n,2)$ to some subcollection of the walks in $\G$.
The graph $\G$ has diameter $d = \lfloor \frac{n}{2} \rfloor$; it has girth four  provided
$n>3$ and $\bar{H}(4,2)$ is a complete bipartite graph.
Let $n>4$ and let $\cK_4$ denote the cell complex with 0-cells and 1-cells given by vertices
and edges of $\G$ and $2$-cells given by the 4-cycles (quadrangles) in $\G$ (cf.\ Definition \ref{Def:cKj}).
\begin{lem}
\label{Lem:4cycles}
Assume $n>4$. Each cycle of length four in $\G$ contains 0 or 2 long edges. For any disk
$B$ in $\cK_4$, the boundary of  $B$ contains an even number of long edges.
\end{lem}
\begin{proof}
If $p=(x_0,x_1,x_2,x_3,x_0)$ is a closed walk in $\G$ where $\{x_0,x_1\}$ is a long edge, then, in the $(d-1)$-cube,
$\partial(x_1,x_0)>3$ and at least one of the remaining edges along $p$ is a long edge in $\G$. As no two long edges
share an endpoint, the number of long edges in any 4-cycle cannot exceed two.
Now consider a disk, composed of 4-cycles and, inductively verify that the mod two sum of the number
of long edges on the boundary remains zero modulo 2 as the number of quadrangles increases.
\end{proof}

\begin{thm}
\label{Thm:foldcube}
Let $n>4$,  let $\G=(X,R)$ be the folded $n$-cube  $\bar{H}(n,2)$, and let $x\in X$.
Consider the corresponding cell complex $\cK_j$, $4\le j\le n-1$.
Then $\cK_{n-1}$ is simply connected while
$\pi\left(\cK_j,x\right) \cong \ints/2\ints$ for $4\le j < n-1$. \qed
\end{thm}

\begin{proof}
As with the folded Johnson graph we will define a closed walk $\bar{p}$ in $\G$ that is the image under $\kappa$ from (\ref{Eqn:CubeCover}) of a shortest path in $H(n,2)$ from one vertex to its antipode.
For $0\le i\le n$, define
$x_i=\underbrace{1\cdots 1}_{i} \, 0\cdots 0$. Then
$$\bar{p}=( \kappa(x_0),\kappa(x_1),\ldots,\kappa(x_{n-1}), \kappa(x_n))$$
is a closed walk in $\G$ of length $n$ based at $x=\kappa(x_0)$ containing exactly one long edge.
 Since any disk composed of $2$-cells in the cell
complex $\cK_4$ coming from $\G$ contains an even number of long edges along its boundary,
this cycle  $\bar{p}$ is not contractible in $\cK_4$. So $[\bar{p}] \neq e$ in $\pi(\cK_4,x)$.

Next, suppose $4<j<n-1$. Theorem \ref{Thm:Hdq} tells us that any $2$-cell in $\cK_j$ has boundary
equal to the boundary of some disk composed entirely of quadrangles. So $[\bar{p}] \neq e$ in
$\pi(\cK_j,x)$ for $j<n$.

On the other hand, if $p$ and $q$ are any closed walks in $\G$ based at $x$ such that $[p]\neq e$ and $[q]\neq e$, we will now show that $[pq] = e$ in $\pi(\cK_4,x)$. It  then follows that $\pi(\cK_4,x) \cong \ints/2\ints$ whenever $n\ge 5$. That is, for $n \ge 5$, $\pi(\G,x,4)$ is a subgroup  of index two in $\pi(\G,x)$.

Let $p$ be any closed walk in $\G$ based at $x$. Decompose $p$ as the concatenation of walks
\begin{equation}
\label{Eqn:Hpqqq}
p = q_1 q_2 \cdots q_s
\end{equation}
where each $q_j$ is a closed walk based at $x$ and does not include $x$ as an internal vertex. If $[q_j] \neq e$, we may lift any $q_j$ under the map (\ref{Eqn:CubeCover}) to a path in $H(n,2)$ in two ways: either from $x_0$ to its antipode $x_n$ or from $x_n$ to $x_0$; otherwise $q_j$ is contractible in $\cK_4$.  Now suppose, for $p,q\in \psi(\G,x)$,
that $p$ and $q$ are not contractible in $\cK_4$ and neither contains $x$ as an internal vertex. Then
the concatenation $pq$ can be viewed as a projection $pq= \kappa(p'q')$ of a closed walk in $H(n,2)$
based at $x_0$ which begins with a walk $p'$ from $x_0$ to its antipode $x_n$ satisfying
$\kappa(p')=p$ followed by a walk $q'$ from $x_n$ to $x_0$ satisfying $\kappa(q')=q$.
Now we use Theorem \ref{Thm:Hdq}
to conclude that $p'q'$ is null homotopic in the cell complex $\cK\left( H(n,2), 4\right)$. Since
$d>4$, each $2$-cell in $\cK\left( H(d,2), 4\right)$ projects, via $\kappa$, to a $2$-cell in $\cK_4$, we
have that $[pq]=e$ in $\pi( \cK_4,x)$. The rest follows from Lewis's Theorem \ref{Thm:Lewis}{\it (ii)} since
the covering diameter of the folded $n$-cube is equal to $n$.
That is,  in the cell complex $\cK_n$, each $q_j$ in (\ref{Eqn:Hpqqq})
is null homotopic giving us that $\pi( \cK_n,x)$ is the trivial group.
\end{proof}

\begin{remark}
Note that, when $n$ is odd, the folded $n$-cube contains cycles of length $n$ whereas
any disk in $\cK_j$, $j<n$,  has a boundary with an even number of edges. For $n$ even,
we needed another approach to show that $G_j \neq \pi(\G,x)$.
\end{remark}

\subsection{Halved cubes}
\label{Subsec:halvedcube}

The halved $n$-cube $\frac{1}{2}Q_n$ is defined as follows. The vertex set $X$ of $\frac{1}{2}Q_n$ is the subset of vertices of the $n$-cube $Q_n = H(n,2)$ (which we think of as ordered binary $n$-tuples) with even Hamming weight. Two vertices $x$ and $y$ are adjacent in $\frac{1}{2}Q_n$ if and only if they are at distance $2$ in $H(n,2)$.

We aim to prove the following theorem.

\begin{thm}
\label{Thm:halfQ}
Let $n\ge 3$, and let $\G=\frac{1}{2}Q_n$ be the halved $n$-cube. Then, for any basepoint $x$, $
\pi(\G,x,3)=\pi(\G,x)$. Equivalently, the fundamental group $\pi(\G,x)$ is generated by equivalence classes of closed walks of essential length at most three. Consequently, the cell complex
$
\cK_3=\cK(\G,3)
$
is simply connected.
\end{thm}

\begin{proof}
The halved cube is vertex transitive, so we may assume that $x=00\cdots 0$. Let $p=(x=x_0,x_1,\dots,x_l=x)$ be a closed walk based at $x$. As before, we will use induction on $c(p)$. If $l(p)\le 3$ then $[p]\in \pi(\G,x,3)$, so we may assume that $l(p)>3$ and $c(p)>0$. Let $\ell=\max\{\partial(x,x_i):x_i\text{ in }p\}$ and $m=\min\{i:\partial(x,x_i)=\ell\}$, and set $a=x_{m-1}, b=x_m, c=x_{m+1}$.

\bigskip

\noindent {\bf Case 0:} $a=c$. \\
$p$ contains the subwalk $(a,b,a)$ which may be removed to create a closed walk $p'=(x=x_0,x_1,\ldots,x_{m-1},x_{m+2},\ldots,x_l=x)$ which is homotopic to $p$ with $c(p')<c(p)$.

\bigskip

\noindent {\bf Case 1:} $a$ and $c$ are adjacent. \\
$(a,b,c,a)$ is a triangle in $\G$. Let
$r=(x_0,x_1,\dots,x_{m-1})$ and $s=(x_{m+1},x_{m+2},\dots, x_l)$. Define
$
q_0 = r(a,b,c,a)r^{-1}
$
and
$
p_0 = r(a,c)s.
$
Then $q_0$ has essential length three, $c(p_0)<c(p)$, and $p\sim q_0p_0$.

\bigskip

\noindent {\bf Case 2:} $a$ and $c$ are not adjacent. \\
In this case, $\partial(a,c)=2$, so $a$ and $c$ differ in exactly four places. Without loss of generality, we may assume
\begin{eqnarray*}
  b &=& \underbrace{1\cdots 1}_{\ell-2} \,110\cdots 0 \\
  a &=& \underbrace{1 \cdots 1}_{\ell-2} \, 000\cdots 0 ~.
\end{eqnarray*}
There are two possibilities for the vertex $c$.

\bigskip

\noindent {\bf Case 2.a:} $c= \underbrace{1\cdots 00 \cdots 1}_{\ell-2} \,110\cdots 0$. \\
Suppose that $c$ the first two zeros in $c$ appear at indices $i$ and $i+1\le \ell-2$. Then $a$, $b$, and $c$ differ only at indices $i,i+1,\ell-1,$ and $\ell$, and are exactly the same elsewhere. To simplify
our notation, we will omit the indices at which all three vertices agree, so we now write $a=1100$, $b=1111$, and $c=0011$. Note that the induced subgraph on the common neighbors of $a$ and $b$ form an
octahedron. The vertices common to the neighborhoods of $a$, $b$, and $c$ are $1010$, $0101$, $1001$, and $0110$. We let $y$ be any one of these four vertices, and notice that $\partial(x,y)<\partial(x,b)$. Let $r=(x_0,x_1,\dots, x_{m-1})$ and $s=(x_{m+1},x_{m+2},\dots, x_l)$, and define the walks
\begin{eqnarray*}
    q_0 &=& r(a,b,y,a)r^{-1} \\
    q_1 &=& r(a,y,b,c,y,a)r^{-1} \\
    p_0 &=& r(a,y,c)s ~ .
\end{eqnarray*}
The walks $q_0$ and $q_1$ each have essential length three, and
\begin{align*}
q_0q_1p_0 &= r(a,b,y,a)r^{-1}r(a,y,b,c,y,a)r^{-1}r(a,y,c)s \\
    &\sim r(a,b,y,a,y,b,c,y,a,y,c)s \\
    &\sim r(a,b,c)s = p.
\end{align*}
Since $\partial(c,y)<\partial(x,b)$, we have $c(p_0)<c_(p)$.

\bigskip

\noindent {\bf Case 2.b:} $c=\underbrace{1\cdots 0 \cdots 1}_{\ell-2} \, 110 \cdots 1 \cdots 0$. \\
Suppose the first 0 entry of $c$ appears at index $i\le \ell-2$ and later has a one at some
index $j\ge \ell+1$. Tuples $a$, $b$, and $c$ are identical except at indices
$i$, $\ell-1$, $\ell$, and $j$, and so we project onto these four relevant coordinate positions and
write $a=1000$, $b=1110$, and $c=0111$. The common neighbors of all three are $0100$, $0010$, $1101$, and $1011$. The vertices denoted by $0100$ and $0010$ are closer to $x$ than $b$ is, so we let
$y$ be either $0100$ or $0010$. We construct $q_0$, $q_1$, and $p_0$ as in Case 2.a, and once again find $p\sim q_0q_1p_0$ where $q_0$ and $q_1$ each has essential length three and $c(p_0)<c(p)$.
\end{proof}

\section{Antipodality and pseudoquotients}
\label{Sec:pseudo}

We have seen a contrast between Hamming and Johnson graphs, on the one hand, and the folded
graphs in cases  where a particular Hamming or Johnson graph is antipodal. In \cite{lewis}, this was
anticipated. While Lewis proved that any $Q$-polynomial distance-regular graph other than a polygon
has girth at most six, she established that  the (group elements represented by) closed walks of essential
length six or less generate the entire fundamental  group except possibly in the case where $\G$ is
a pseudoquotient. We have exhibited two families of quotients  where this exceptional behavior occurs.
Here we briefly recap the portion of Lewis's paper dealing with pseudoquotients. As we are not using
the parameters of the Bose-Mesner algebra or $Q$-sequences of $Q$-polynomial distance-regular
graphs, our coverage will not be self-contained but will rather refer to material from \cite{lewis} and \cite{bcn}.

A distance-regular graph $\G$ of diameter $D$ is \emph{antipodal} if the relation $x\sim y$ if
$\partial(x,y)\in \{0,D\}$ is an  equivalence relation. Examples include hypercubes $H(D,2)$,
the Johnson graphs $J(2D,D)$, and polygons on $2D$ vertices, as well as the icosahedron,
dodecahedron and the line graph of the Petersen graph. Suppose now that $\G$ is an antipodal
distance-regular graph of diameter $D\ge 3$. Identifying antipodal vertices and joining two antipodal
classes by an edge whenever one contains a neighbor of some vertex in the other, we  obtain
another distance-regular graph called the \emph{folded graph} and often denoted $\bar{\G}$.
In \cite[Section 4.2]{bcn}, we see that the spectral information about $\bar{\G}$ is closely related to
that of $\G$. Lewis's concept of a  pseudoquotient is motivated by the relationship between the
$Q$-sequences of these two graphs when both are $Q$-polynomial.

Let $\G =(X,R)$ be a distance-regular graph of diameter $d$  with adjacency matrix $A$. For $0\le i\le d$, define the ``distance $i$'' matrix $A_i \in \re^{X\times X}$  of $\G$ via $(A_i)_{xy}=1$ if $\partial(x,y)=i$ and  $(A_i)_{xy}=0$ otherwise.  (So $A_0=I$, $A_1=A$, $\sum_i A_i=J$, etc.) Assume matrix $A$ has  spectral decomposition
\begin{equation}
\label{Eqn:spec}
A = \sum_{\theta} \theta E_\theta
\end{equation}
 where $E_\theta^2 =  E_\theta$, $E_\theta E_\tau = 0$ for $\theta\neq  \tau$, and the eigenvalues in the above sum are pairwise distinct \cite[Section 4.1B]{bcn}.
 The distance-regular graph  $\G$ is $Q$-\emph{polynomial} if one can find an ordering $E_0,E_1,\ldots,E_d$ of the idempotents
 $\{E_\theta : \theta\}$ such that
 \[
 \left\langle E_1^{\circ 0}, E_1^{\circ 1},\ldots,E_1^{\circ j} \right\rangle =
  \left\langle E_{0}, E_1,\ldots,E_{j} \right\rangle
  \]
  for every $j$, $0\le j\le d$ where $E^{\circ n}$ is the matrix obtained by raising each entry of $E$ to the $n^{\rm th}$ power. Clearly
  this is equivalent to the condition that, for each $j$, $E_j$ is expressible as a polynomial of degree exactly $j$, applied entrywise,
  in the matrix $E_1$. In this case, we call $E_1$ a $Q$-\emph{idempotent}. Note that a $Q$-polynomial distance-regular graph may
  have several $Q$-idempotents.

  Now assume that $E_1$ is a $Q$-idempotent for $\G$. As we have two bases, $\{A_0,\ldots,A_d\}$ and $\{E_0,\ldots,E_d\}$ for the adjacency algebra   of $\G$, there are unique \emph{dual eigenvalues} $\theta_0^*,\ldots,\theta_d^*$ satisfying
  \begin{equation}
\label{Eqn:dualeigs}
E_1 = \sum_{i=0}^d \theta_i^* A_i ~.
\end{equation}
In this case, the sequence $\theta_0^*,\ldots,\theta_d^*$ is called a $Q$-\emph{sequence} for $\G$. Note that, as there
may be multiple $Q$-idempotents, $\G$ may correspondingly have multiple $Q$-sequences. Douglas Leonard proved (cf.\ \cite[Lemma 31]{lewis}) that, for any $Q$-sequence $\theta_0^*,\ldots,\theta_d^*$, $d\ge 3$, there is a unique nonzero $\lambda$ satisfying
  \begin{equation}
\label{Eqn:dougslambda}
\theta_{i-2}^* - \theta_{i-1}^* = \lambda \left(\theta_{i-3}^* - \theta_{i}^* \right), \qquad 3\le i\le d
\end{equation}

If $\G$ happens to be the quotient of an antipodal distance-regular graph of diameter $D\ge 7$, then $\G$ has
diameter $d=\lfloor D/2 \rfloor$ and every $Q$-sequence
for $\G$ may be extended via the rule $\theta_{i}^*=\theta_{D-i}^*$ ($d+1\le i\le D$) in such a way that (\ref{Eqn:dougslambda}) still
holds for these additional terms in the sequence.

\begin{defn}
 \label{Def:pseuodoquotient}
 Let $\G$ be a $Q$-polynomial distance-regular graph of diameter $d\ge 3$. Then $\G$ is a \emph{pseudoquotient}
 if there is an integer $D\in \{2d,2d+1\}$ such that, for every $Q$-sequence $\theta_0^*,\ldots,\theta_d^*$ of $\G$,
the extended sequence given by the rule  $\theta_{i}^*=\theta_{D-i}^*$ ($d+1\le i\le D$) satisfies
\[
\theta_{i-2}^* - \theta_{i-1}^* = \lambda \left(\theta_{i-3}^* - \theta_{i}^* \right), \qquad 3\le i\le D~.
\]
\end{defn}

\begin{example}
For example, a known open parameter set for a potential distance-regular graph of diameter three has intersection array
$\{49,40,22;1,5,28\}$. If such a graph exists, it has $Q$-sequence $112, 32, 2, -8$ satisfying (\ref{Eqn:dougslambda}) with
$\lambda = 1/4$. In fact, the sequence $112, 32, 2, -8, -8, 2, 32, 112$ also satisfies  (\ref{Eqn:dougslambda}) with $\lambda=1/4$
but there can be no distance-regular graph with intersection array $\{ 49,40,22, t ,28,5,1; 1,5,28, t, 22,40,49\}$ (where $t$ is
some integer) for the obvious reason that $p_{1,3}^2 = 22 < p_{1,5}^4 = 28$, violating an elementary feasibility condition.
\end{example}

\section{Further directions}
\label{Sec:further}
The Johnson graph and Folded Johnson graphs are locally grid graphs. One might investigate other locally grid graphs even if they are not distance regular.

Theorem~\ref{Thm:Jnk} can also be used to show that a Johnson graph
of diameter at least two cannot be obtained by identifying the
antipodal vertices of a larger distance-regular graph. Can our
decompositions of closed walks similarly be used to rule out
antipodal covers of other distance-regular graphs, or to show
that their known antipodal covers are unique?

There are several families of graphs closely related to the Johnson and Hamming graph that are worthy of study in this topological context. The \emph{Odd graph} $O_k$ is the distance-($k-1$) graph of the Johnson graph $J(2k-1,k-1)$ and this
graph played a role in our analysis of the folded Johnson graphs. One may also consider the bipartite double of $O_k$,
the \emph{doubled Odd graph}, or distance-regular bipartite doubles in general. In the case $k=3$, and only in this case,
there is a second antipodal cover: in addition to the Desargues graph (the bipartite double of the Petersen graph $O_3$),
we also have the dodecahedron, which is not $Q$-polynomial. Can we rule out non-distance regular antipodal covers of Odd graphs using homotopy?(Cf. \cite[Section 9]{godsilhensel}).

With one trivial exception, the Hamming graphs have no distance-regular antipodal covers. (The $8$-cycle is an antipodal double cover of $H(2,2)$.)  The folded $4$-cube has one non-Hamming antipodal cover as does the folded $5$-cube. In all other cases the only antipodal covers of folded cubes are cubes themselves. We have not computed $g_2(\G)$ for either of these graphs. The overall goal is to
better understand the parameter $g_2(\G)$, especially in the non-$Q$-polynomial case, and its relationship
to conjectures about the girth of distance-regular graphs.

\section*{Acknowledgments}

The authors thank Herman Servatius for helpful conversations on material in this paper. The second author’s earlier work on this topic began with support from the National Science Foundation which is gratefully acknowledged.

\end{document}